\documentclass[11pt]{amsart}
\usepackage{fullpage,latexsym,amsmath,amssymb,amsfonts,amscd,graphics}
\usepackage[mathscr]{eucal}
\usepackage[all]{xypic}
\numberwithin{equation}{section}
\theoremstyle{plain}
\newtheorem{theorem}[equation]{Theorem}

\newtheorem{lemma}[equation]{Lemma}
\newtheorem{proposition}[equation]{Proposition}
\newtheorem{corollary}[equation]{Corollary}

\theoremstyle{remark}
\newtheorem{remark}[equation]{Remark}

\theoremstyle{definition}
\newtheorem{definition}[equation]{Definition}
\newtheorem{notation}[equation]{Notation}

\newcommand{\bP}{\mathbb{P}}
\newcommand{\bH}{\mathbb{H}}

\newcommand{\bQ}{\mathbb{Q}}
\newcommand{\bZ}{\mathbb{Z}}
\newcommand{\bF}{\mathbb{F}}
\newcommand{\bC}{\mathbb{C}}
\newcommand{\bL}{\mathbb{L}}

\newcommand{\sH}{\mathrm{H}}

\newcommand{\calF}{\mathcal{F}}

\newcommand{\calH}{\mathcal{H}}
\newcommand{\calN}{\mathcal{N}}
\newcommand{\calT}{\mathcal{T}}
\newcommand{\calh}{\mathfrak{h}}
\newcommand{\fN}{\mathfrak{N}}
\newcommand{\calM}{\mathcal{M}}
\newcommand{\calO}{\mathcal{O}}
\newcommand{\calL}{\mathcal{L}}
\newcommand{\calP}{\mathcal{P}}

\newcommand{\calX}{\mathcal{X}}

\newcommand{\calD}{\mathcal{D}}

\newcommand{\Sym}{\mathrm{Sym}}

\newcommand{\Hom}{\mathrm{Hom}}
\newcommand{\SL}{\mathrm{SL}}
\newcommand{\Proj}{\mathrm{Proj}}
\newcommand{\Sing}{\mathrm{Sing}}

\newcommand{\Gr}{\mathrm{Gr}}

\newcommand{\corank}{\mathrm{corank}}

\newcommand{\gquot}{/\!\!/}

\author{Radu Laza}
\title[Moduli of cubic $4$-folds]{The moduli space of cubic fourfolds via the period map} 
\address{University of Michigan \\
1832 East Hall \\
Ann Arbor, MI 48109}
\email{rlaza@umich.edu}

\begin{document}
\bibliographystyle{amsplain}

\begin{abstract}
We characterize the image of the period map for cubic fourfolds with at worst simple singularities as the complement of an arrangement of hyperplanes in the period space. It follows then that the GIT compactification of the moduli space of cubic fourfolds is isomorphic to the Looijenga compactification associated to this arrangement. This paper builds on and is a natural continuation of our previous work on the GIT compactification of the moduli space of cubic fourfolds.
\end{abstract}

\maketitle
%%%%%%%%%%%%%%%%%%%%%%%%%%%%%%%%%%%%%%%%%%%%%%%%%%%%%%%%%%%%%% 
\section{Introduction}\label{sectintro}
An indispensable tool for the study of  K3 surfaces is the period map. The period map for K3 surfaces  is both injective (satisfies the global Torelli theorem) and surjective. It follows that the moduli space of algebraic K3 surfaces of a given degree is the quotient of a $19$-dimensional bounded symmetric domain by an arithmetic group. This explicit description of the moduli space  has numerous geometric applications, for example to the study of deformations of certain classes of surface singularities (e.g. the work of Pinkham \cite{pinkhamduality}, Looijenga \cite{looijengatriangle}, and others). The study of deformations of higher-dimensional singularities led us to consider the period map for  cubic fourfolds, which is well known to behave similarly to the period map for K3 surfaces. Specifically, as in the case of K3 surfaces, the period domain $\calD$ for cubic fourfolds is a bounded symmetric domain of type IV of the right dimension, and the global Torelli theorem holds  (Voisin \cite{voisin}). In this paper we analyze the remaining open question, the question of surjectivity (or more precisely of characterizing the image) of the period map for cubic fourfolds. A surjectivity type statement is needed in applications, especially the applications to the study of  degenerations and singularities of cubic fourfolds.

\smallskip

Our main result is a characterization of the image of the period map for cubic fourfolds, giving a positive answer to a conjecture of Hassett \cite[\S4.4]{hassett}.  This builds on earlier results of Voisin \cite{voisin} and Hassett  \cite{hassett}, and was inspired by the  recent work of Allcock--Carlson--Toledo \cite{allcock1,allcock3fold} and Looijenga--Swierstra \cite{looijengaswierstra} on the moduli space of cubic threefolds.

\begin{theorem}\label{mainthm1}
The image of the period map for cubic fourfolds $\calP_0:\calM_0\to \calD/\Gamma$ is the complement of the hyperplane arrangement $\calH_{\infty}\cup \calH_\Delta$ (see Def. \ref{defdethyp}). Furthermore, the period map extends to a regular morphism $\calP:\calM\to \calD/\Gamma$ over the simple singularities locus with image the complement of the arrangement $\calH_{\infty}$.
\end{theorem}

We note that the theorem is analogous to the corresponding statement for degree two K3 surfaces. Specifically, the hyperplane arrangement $\calH_{\infty}\cup \calH_\Delta$ is  the analogue of the hyperplane arrangement corresponding to the $(-2)$-curves in the K3 case. 
As in the case of low degree K3 surfaces, the two components $\calH_{\Delta}$ and $\calH_{\infty}$ are distinguished by certain arithmetic properties and they parametrize two kinds of degenerations of cubic fourfolds. Finally, again completely analogous to the  K3 situation (see \cite[\S8]{looijengacompact}),  a stronger result holds: the compactification $\overline{\calM}$ of the moduli of cubic fourfolds obtained by  means of geometric invariant theory (GIT) is an explicit birational modification of the Baily--Borel compactification $(\calD/\Gamma)^*$ of the period space. 
\begin{theorem}\label{mainthm2}
The period map for cubic fourfolds induces an isomorphism 
$$\overline{\calM}\cong \widetilde{\calD/\Gamma},$$ 
where $\overline{\calM}$ is the GIT compactification of moduli space of cubic fourfolds, and $\widetilde{\calD/\Gamma}$ denotes the Looijenga compactification associated to the arrangement of hyperplanes $\calH_{\infty}$.
\end{theorem}

Theorem \ref{mainthm2} follows immediately from Theorem \ref{mainthm1} and the general results of Looijenga (esp. \cite[Thm. 7.6]{looijengacompact}). Thus, we are mainly concerned here with establishing Theorem \ref{mainthm1}. For this we take an incremental approach, following  the arguments of Shah \cite{shahinsignificant,shah} on the surjectivity of the period map for low degree K3 surfaces. We start by computing the GIT compactification  $\overline{\calM}$ of the moduli space of cubic fourfolds. By studying the monodromy around cubic fourfolds with closed orbits, we conclude that the indeterminacy of the period map is a curve $\chi$. Next, we successively blow-up $\overline{\calM}$ first in a special point $\omega\in \chi$, and then in the strict transform of $\chi$. A new monodromy analysis for the blow-up, allows us to conclude that Theorem \ref{mainthm1} holds. Further details on the organization of the paper and the main intermediary results are given below.

\smallskip

The GIT analysis for cubic fourfolds was done in Laza \cite{gitcubic} (some partial results were also obtained by Allcock and Yokoyama). We recall the relevant details in section \ref{sectpreliminary}. Basically, what we need from the GIT computation are the following results. First, a cubic fourfold having at worst simple isolated singularities is GIT stable. In particular, we can talk about the moduli space $\calM$ of such cubic fourfolds. For monodromy reasons, over $\calM$ the period map naturally extends. The second GIT result used is that the boundary of $\calM$ in the GIT compactification $\overline{\calM}$ is  naturally  stratified in 3 types, labeled II, III, and IV. This stratification is closely related to the stratification of  Shah \cite[Thm. 2.4]{shah}. Essentially, the singularities that occur for Type II and III fourfolds are the insignificant limit singularities of Shah \cite{shahinsignificant}. Consequently, the Type II and III fourfolds  cause no problem for the period map. Finally, the GIT results identify the locus of Type IV fourfolds in $\overline{\calM}$ (the indeterminacy locus of the period map) to be a rational curve $\chi$ containing a special point $\omega$. In our analysis, the Type IV fourfolds play the same role as the triple conic for plane sextics (see \cite{shah}) or the chordal cubic for cubic threefolds (see \cite{allcock3fold,looijengaswierstra}). 

\smallskip

Based on the GIT results mentioned above, the proof of the Theorem \ref{mainthm1} follows in two main steps. The first step, completed in section \ref{sectmonodromy}, is to prove that we can control the monodromy of  $1$-parameter degenerations of cubic fourfolds with central fiber not of Type IV (strictly speaking, this is weaker than an extension statement for the period map $\calP$ for Type II and III fourfolds, but since it suffices for our purposes, will only talk about the Type IV locus as the indeterminacy for $\calP$).  There are two  ingredients going in the proof of the previous statement: a reduction to the central fiber $X_0$, followed by a Hodge theoretical computation for $X_0$. Specifically, in \S\ref{reductioncentral}, we prove that the natural specialization morphism associated to a degeneration is injective on certain pieces of the corresponding mixed Hodge structures. Thus,  the question about the monodromy of the family can be reduced to checking some statement about the mixed Hodge structure of the central fiber $X_0$. This in turn is relatively easy in our situation. Namely, the mixed Hodge structure of a singular cubic fourfold $X_0$ can be computed by using the projection from a singular point (see \S\ref{sectmhs}). This reduces the computation  to standard facts about degenerations of $K3$ surfaces. As hinted in the previous paragraph, the essential fact that makes the proofs of section \ref{sectmonodromy} work  is that the singularities of Type II and III fourfolds are double suspensions of special surface singularities, the so-called insignificant limit singularities of Mumford and Shah. 

\smallskip

In section \ref{sectresolution}, we analyze the degenerations to Type IV fourfolds, completing the second step of our proof. A partial analysis of the degenerations to Type IV fourfolds was done by Hassett \cite[\S4.4]{hassett} and Allcock-Carlson-Toledo \cite[\S5]{allcock3fold}. We complete the analysis by proceeding as follows. First, we note  that the indeterminacy locus of the period map, the curve $\chi$   parametrizing the Type IV fourfolds, is the locus where the GIT quotient $\overline{\calM}$ has the worst singularities. Specifically, the Type IV fourfolds are characterized among the semi-stable cubic fourfolds by the fact that their stabilizer is not virtually abelian. Therefore, it is natural to consider a partial desingularization $\widetilde{\calM}$ of $\overline{\calM}$ as constructed by Kirwan \cite{kirwan}. Namely, we let $\widetilde{\calM}$ be the blow-up the special point $\omega\in \chi$, followed by the blow-up of the strict transform of $\chi$  (see \S\ref{sectkirwan}).  The space $\widetilde{\calM}$ has only toric singularities and the  period map essentially extends over $\widetilde{\calM}$. More precisely,  we note that the effect of these two blow-ups is to replace the Type IV fourfolds by some fourfolds, that we call of Type I'--III', having the same type of singularities as the cubic fourfolds of Type I-III. Thus, we can control the monodromy as in  section \ref{sectmonodromy}. We then show that the limit period point corresponding to a Type I' fourfold belongs to $\calH_{\infty}/\Gamma$. This is enough to complete the proof of theorem \ref{mainthm1}.

\smallskip

Finally, in section \ref{sectproofs}, we put everything together and prove the two main theorems \ref{mainthm1} and \ref{mainthm2}. As a simple application, we note that from theorem \ref{mainthm2} and the results of Looijenga \cite{looijengacompact} we can recover some information about the GIT compactification purely in arithmetic terms. This is discussed in the last section, section \ref{sectapplication}.

\smallskip

The study of the periods of cubic fourfolds is closely related to the study of the periods of irreducible holomorphic symplectic fourfolds. Namely, the Fano variety of a cubic fourfold $X$ is a symplectic fourfold $F$ deformation equivalent to the resolution of the symmetric square of a $K3$ surface (see Beauville--Donagi \cite{beauvilledonagi}). The periods of $X$ are essentially the periods of $F$. It follows that the period domain $\calD/\Gamma$ can be interpreted also as the period domain of such symplectic fourfolds (with a degree $6$ polarization). A surjectivity type statement for irreducible symplectic fourfolds is known (see Huybrechts \cite{huybrechts}), but the linear systems on symplectic fourfolds are not enough understood to obtain a characterization  of the image of the period map in this way. We note however that, from the perspective of symplectic fourfolds, the complement of the image of the period map, the divisor $\calH_\infty/\Gamma$, has a geometric meaning: it is the locus of symmetric squares of degree two $K3$ surfaces in the period domain. These type of fourfolds are not Fano varieties of cubic fourfolds, but rather of certain quadric bundles associated to   degree two K3 surfaces   (see \S\ref{sectmodular}). Thus, going back to cubics, we can interpret the Kirwan blow-up  as enlarging the moduli space of cubic fourfolds $\overline{\calM}$ to include these quadric bundles. The first blow-up corresponds to general degree two K3 surface, and the second to the special (i.e. unigonal) ones.

\smallskip

%The Fano variety of lines on a cubic fourfold $X$ is an irreducible symplectic fourfold $F$ deformation equivalent to the $\Hilb^2(S)$ for $S$ a K3 surface (\cite{beauvilledonagi}). The periods of $X$ are essentially the periods of $F$. For irreducible symplectic fourfolds the surjectivity of the  period map is known. However, the linear systems on symplectic fourfolds are not well enough understood (see \cite{hassetttschinkel}) to make possible a characterization of the image of the period map for cubic fourfolds starting from the Fano variety. On the other hand we note that locus $\calH_{\infty}/\Gamma$ excluded from the image of the period map for cubic fourfolds corresponds precisely to irreducible symplectic fourfolds of  type $\Hilb^2(S)$ for $S$ a degree two K3 surface. These type of fourfolds are not Fano varieties of cubic fourfolds, but rather of certain quadric bundles associated to   degree two K3 surfaces   (see \S\ref{sectmodular}). The Kirwan blow-up enlarges the moduli space of cubic fourfolds $\overline{\calM}$ to include these quadric bundles. The first blow-up corresponds to generic degree two K3 surface, and the second to the special (i.e. unigonal or elliptic) ones (see \S\ref{sectkirwan}). In conclusion, we essentially obtain surjectivity for the period map for cubic fourfolds compatible with the surjectivity for irreducible symplectic fourfolds.

%\smallskip

An alternative independent proof of our main result was obtained by Eduard Looijenga \cite{looijengacubic}. While the proof of \cite{looijengacubic} also uses some of the GIT results, the techniques (based on \cite{looijengaswierstra1}) to handle the period map for degenerations of cubic fourfolds are different. We would like to thank Eduard Looijenga for informing us about his work. 

\subsection*{Notations and Conventions} Our notations and terminology are based on Mumford \cite{GIT} when we refer to GIT, on Arnold et. al. \cite{AGV1} when we refer to singularities, and on Griffiths et al. \cite{griffithsbook} when we refer to Hodge Theory. Additionally, we are using freely the notations of Laza \cite{gitcubic} (esp. \cite[\S1.2]{gitcubic}). In particular, $\calM_0$, $\calM$, and $\overline{\calM}$ denote the moduli space of smooth cubic fourfolds, of cubic fourfolds with simple (A-D-E) singularities, and the GIT compactification respectively. 

All hypersurface singularities are considered up to stable equivalence (i.e. up to adding squares of new variables to the defining equation). In particular, it makes sense to say that  two hypersurface singularities of different dimension have the same analytic type.   

%%%%%%%%%%%%%%%%%%%%%%%%%%%%%%%%%%%%%%%%%%%%%%%%%%%%%%%%%%%%%% 
\section{Preliminary results}\label{sectpreliminary}
In this section we collect a series of results that we will need in the subsequent sections. Some are well known general results, others are specific to cubic fourfolds and are based mostly on Voisin \cite{voisin}, Hassett \cite{hassett}, and Laza \cite{gitcubic}. 

\subsection{The GIT compactification of the moduli space of cubic fourfolds}\label{gitresults}  
The computation of the GIT compactification $\overline{\calM}$ of the moduli space of smooth cubic fourfolds $\calM_0$ was carried out in Laza \cite{gitcubic}. Here we recall the GIT results that are needed in what follows (for details see \cite{gitcubic}).
 
\smallskip

We start by noting that a cubic fourfold with only simple singularities is stable.  

\begin{theorem}
A cubic fourfold with at worst simple  singularities is GIT  stable. In particular, it makes sense to talk about the moduli space $\calM$ of such fourfolds as an open subset of the GIT quotient $\overline{\calM}$. The boundary $\overline{\calM}\setminus \calM$ consists of six closed irreducible components, that we label $\alpha$--$\phi$. The general point of a boundary component corresponds to a semistable cubic fourfold (with closed orbit) with singularities as given in table \ref{tableboundary}.
\end{theorem}

\begin{table}[htb]
\begin{center}
\renewcommand{\arraystretch}{1.25}
\begin{tabular}[2cm]{|c|c|l|}
\hline
Dim.&Component & Singularities of the cubics parameterized  by the boundary component \\
\hline\hline
1&$\alpha$& an elliptic normal curve of deg. 4 and a rational normal curve of deg. 1\\
\hline
3&$\beta$& two isolated $\widetilde{E}_8$ singularities\\
\hline
2&$\gamma$& one isolated $\widetilde{E}_7$ singularity and an elliptic normal curve of deg. 2\\
\hline
1&$\delta$&three isolated $\widetilde{E}_6$ singularities\\
\hline
3&$\epsilon$&a  rational normal curve of degree 4\\
\hline
2&$\phi$& an elliptic normal curve of degree 6\\
\hline
\end{tabular}
\vspace{0.1cm}
\caption{The boundary  of $\calM$ in $\overline{\calM}$}\label{tableboundary}
\end{center}
\end{table}

Explicit equations for the cubic fourfolds with closed orbit parameterized  by the boundary components 
are given in \cite{gitcubic}. Here we only need that, with only one exception (the Type IV described below), the singularities occurring for these cubics are the insignificant limit singularities of Mumford and Shah. Specifically, we recall the following list of singularities from Shah \cite[Thm. 1] {shahinsignificant}:

\begin{definition}\label{shahlist}
We say that a $2$-dimensional hypersurface singularity $(X,0)\subset (\bC^3,0)$ is of type $t_1$--$t_6$ if its defining equation is one of the following: 
\begin{itemize}
\item[($t_1$)]: the equation of an isolated rational double point, i.e.  $A_n$, $D_m$, or $E_r$;
\item[($t_2$)]: $x_2x_3$, i.e. double line or $A_\infty$;
\item[($t_3$)]: $x_1^2x_2+x_3^2$, i.e. ordinary pinch point or $D_\infty$;
\item[($t_4$)]: $x_3^2+(x_1+a_1x_2^2)(x_1+a_2x_2^2)(x_1+a_3x_2^2)+g(x_1,x_2)$ with $\mathrm{ord}_0( g )>6$ with respect to the weights $2$ and $1$  for $x_1$  and $x_2$ respectively, and such that at least two of the $a_i$ are distinct;
\item[($t_5$)]: $x_3^2+f_4(x_1,x_2)+g(x_1,x_2)$ with $f_4$ a homogeneous polynomial of degree $4$ and $\mathrm{ord}_0( g )>4$, and such that   $f_4$ has no triple root;
\item[($t_6$)]: $f_3(x_1,x_2,x_3)+g(x_1,x_2,x_3)$ with $f_3$ a homogeneous polynomial of degree $3$ and $\mathrm{ord}_0( g )>4$, and such that $f_3$ has at worst ordinary double points as singularities.
\end{itemize}
where $g$ is a convergent power series in the appropriate variables.
\end{definition}

\begin{remark}\label{reminsignficant} 
The types $t_4$--$t_6$ cover the simple elliptic (denoted $\widetilde{E_r}$), cusp, and degenerate cusp singularities. In fact, the singularities occurring in the list are the $2$-dimensional  semi-log-canonical hypersurface singularities (cf. \cite[Ch. 15]{AGV1}  and \cite[Thm. 4.21]{ksb}). The relevance of this fact is that for $K3$ surfaces three different concepts of measuring the complexity of a singularity are almost equivalent: semi-log-canonical, cohomologically insignificant (cf. \cite[Thm. 1, Thm. 2] {shahinsignificant} and \cite{dolgachevinsignificant,steenbrinkinsignificant}), and GIT stable (see \cite[\S3]{mumford}). This gives a conceptual explanation of the close relation between the GIT and Hodge theoretical constructions of the moduli spaces of low degree K3 surfaces. 
\end{remark}

We adapt the previous definition for fourfold singularities as follows:
\begin{definition}\label{shahlist2}
We say that a $4$-dimensional hypersurface singularity $(X,0)\subset (\bC^5,0)$ is of type $t_1$--$t_6$ if it is a double suspension of a surface singularity of type  $t_1$--$t_6$  respectively (i.e. the defining equation can be taken of the form $f(x_1,x_2,x_3)+x_4^2+x_5^2$ with $f$ of type  $t_1$--$t_6$). 
\end{definition}

In analogy to the work of Shah \cite{shahinsignificant,shah} on degenerations of low degree $K3$ surfaces, it is natural to introduce the following stratification in types of the cubic fourfolds:
\begin{definition}\label{deftypecubic}
Let $X$ be a semistable cubic fourfold. We say that
\begin{itemize}
\item[i)] {\it $X$ is of Type I} if $X$ has at worst simple singularities (N.B. $X$ is stable);
\item[ii)] {\it $X$ is of Type II} if $X$ is not of Type I and all its singularities are of type $A_n$, $D_m$, $E_r$, $A_\infty$ (type $t_2$), $D_\infty$ (type $t_3$), or $\widetilde{E_r}$ for $r=6,\dots 8$ (the generic case of type $t_4$--$t_6$); 
\item[iii)] {\it $X$ is of Type III} if $X$ is not of Type I or II and all its singularities are of type $t_1$--$t_6$;
\item[iv)] {\it $X$ is of Type IV} if $X$ is not of Type I, II, or III. 
\end{itemize}
\end{definition}

By definition the singularities of Type I--III fourfolds are quite mild. Consequently one can control the monodromy of $1$-parameter degenerations to Type I--III fourfolds. In fact, the type is simply the index of nilpotency of the monodromy. In contrast, the degenerations to Type IV fourfolds are more complicated and should be treated separately. Thus, we will need the following GIT result that identifies the locus of Type IV fourfolds:

\begin{theorem}
The semistable cubic fourfolds of Type IV with closed orbits are parameterized  by a rational curve $\chi$ (including the special point $\omega$). Furthermore, the Type IV fourfolds are characterized among the semistable cubics with closed orbits by the fact that their stabilizer is not virtually abelian; the stabilizer is  $\SL(2)$ for a general point on $\chi$ (resp. $\SL(3)$ for $\omega$). 
\end{theorem}

The  Type IV fourfolds parametrized by the curve $\chi$ are given by the equations:
\begin{equation}\label{eqchi}
\chi:\ \ g(x_0,\dots,x_5)=\left| \begin{array}{ccc}
x_0&x_1&x_2\\
x_1&x_2&x_3\\
x_2&x_3&x_4
\end{array}\right|+ax_5(x_0 x_4-4 x_1 x_3 +3 x_2^2)+b x_5^3,
\end{equation}
for $(a:b)\in W\bP(1:3)$. As a special case, one obtains the determinantal cubic (the case labeled $\omega$):
\begin{equation}\label{eqomega}
\omega:\ \ g(x_0,\dots,x_5)=\left| \begin{array}{ccc}
x_0&x_1&x_2\\
x_1&x_5&x_3\\
x_2&x_3&x_4
\end{array}\right|.
\end{equation}
The fourfold of type $\omega$ is the secant to the Veronese surface in $\bP^5$. In particular, it is singular along the Veronese surface. More generally, referring to \cite{gitcubic}, we note the following:

\begin{proposition}\label{propsingtype}
Let $X$ be a semistable cubic fourfold with closed orbit.
\begin{itemize}
\item[i)] If $X$ is of Type II, then the non-simple singularities of $X$ are of the following types:
\begin{itemize}
\item[-] isolated singularities of type $\widetilde{E_r}$; 
\item[-] a rational normal curve $C$ of singularities; at a general point of $C$ the singularity is of type $A_\infty$; at $4$ special points on $C$ the singularities are of type $D_\infty$;
\item[-] an elliptic normal curve $C$ of singularities; at every point of $C$ the singularity is of type $A_{\infty}$.
\end{itemize}
Additionally, the combinations of non-simple singularities that occur are given in table \ref{tableboundary}.
\item[ii)] If $X$ is of Type III, then the locus of non-simple singularities consists of a connected union of rational normal curves. At a general point of a curve of singularities, the singularity  is of type $A_\infty$. On each such curve, there exist two special points with singularities of type $t_4$--$t_6$. Furthermore, the only case when all the curves of singularities are lines is the case $\zeta$ described below.
\item[iii)] If $X$ is of Type IV, then, except the case $\omega$, $X$ is singular along a rational normal curve $C$ of degree $4$ with $A_2$ transversal singularities. If $X$ is the fourfold $\omega$, then $X$ is singular along  a Veronese surface in $\bP^5$. Furthermore, this is the only case of a semistable cubic fourfold with $2$-dimensional singular locus.
\end{itemize}
\end{proposition}

As suggested by the proposition, the most degenerate case of fourfold of Type II or III is the toric fourfold (labeled $\zeta$):
\begin{equation}\label{eqzeta}
\zeta:\ \ g(x_0,\dots,x_5)=x_0x_4x_5+x_1x_2x_3.
\end{equation}
The fourfold of type $\zeta$ is singular along $9$ lines, meeting in triples. We also note that the point $\zeta$ in $\overline{\calM}$ (corresponding to the orbit of the fourfold given by (\ref{eqzeta})) is the intersection of all the boundary strata $\alpha$--$\phi$.

\begin{remark}
Some of the Type IV fourfolds appeared previously in literature in connection to the period map. Namely, the case $\omega$ was considered by Hassett \cite[\S4.4]{hassett}, and the special case of $\chi$ with extra $\mu_3$-automorphisms  (i.e. $a=0$ in (\ref{eqchi})) by Allcock-Carlson-Toledo \cite{allcock3fold}.
\end{remark}

%%%%%%%%%%%%%%%%%%%%%%%%%%%%%%%%%%%%%%%%%%%%%%%%%%%%%%%%%%%%%%
\subsection{The period map for cubic fourfolds}\label{sectperiod}
The period map for cubic fourfolds is defined by sending a smooth cubic fourfold $X$ to its periods:
$$\calP_0:\calM_0\to \calD/\Gamma.$$
The space $\calD$ is the classifying space of polarized Hodge structures on the middle cohomology of a cubic fourfold, and $\Gamma$ is the monodromy group. 

\begin{notation}
We denote  $\Lambda:=\langle1\rangle^{\oplus 21}\oplus \langle-1\rangle^{\oplus 2}$ the abstract lattice isometric to the integral cohomology of a cubic fourfold, by $h\in \Lambda$ the polarization class (the square of the class of a hyperplane section),  and by $\Lambda_0:=\langle h \rangle^\perp\cong E_8^{\oplus 2}\oplus U^{\oplus 2}\oplus A_2$ the primitive cohomology. 
\end{notation}

The Hodge numbers of a cubic fourfold are $h^{4,0}=h^{0,4}=0$, $h^{3,1}=h^{1,3}=1$, and $h^{2,2}=21$. Thus, the classifying space of Hodge structures for cubic fourfolds is a $20$-dimensional Type IV bounded symmetric domain $\calD\cong\textrm{SO}_0(20,2)/\textrm{S}\left(\textrm{O}(20)\times\textrm{O}(2)\right)$.  
We regard $\calD$ as a choice of one of the two connected components of the space of lines in $\Lambda_0\otimes\bC$ satisfying the two Riemann--Hodge bilinear relations:
$$\calD=\{\omega\in\bP(\Lambda_0\otimes_\bZ \bC)\mid \omega^2=0,\ \omega.\bar{\omega}<0\}_0$$ 
(N.B. the subscript $0$ indicates the choice of a connected component). By results of Ebeling and Beauville \cite[Thm. 2]{beauvillemon}, $\Gamma$  can be canonically identified to both $\mathrm{O^+_h}(\Lambda)$ (the automorphisms of $\Lambda$ that preserve the polarization and the orientation of a negative definite $2$-plane in $\Lambda$) and  $\mathrm{O^*}(\Lambda_0)$ (the automorphisms of $\Lambda_0$ that act trivially on the discriminant group $A_{\Lambda_0}:=(\Lambda_0)^*/\Lambda_0$ and preserve the orientation of a negative definite $2$-plane in $\Lambda_0$). Note that $\Gamma\cong \mathrm{O^*}(\Lambda_0)$ respects the choice of connected component for $\calD$, and thus acts naturally on $\calD$.  The orbit space $\calD/\Gamma$ is a quasi-projective variety with a natural compactification, the Baily--Borel compactification $(\calD/\Gamma)^*$.

\smallskip

It is a general fact that the period map for cubic fourfolds is a local isomorphism. In fact, the stronger global Torelli theorem holds in this situation by a result of Voisin \cite{voisin}. Since the period map is algebraic, it follows that  $\calP_0$ is  a birational isomorphism between the quasi-projective varieties $\calM_0$ and $\calD/\Gamma$. For us, it is more convenient to regard $\calP_0$ as a birational map between the projective varieties $\overline{\calM}$  and $(\calD/\Gamma)^*$ that compactify $\calM_0$ and $\calD/\Gamma$ respectively.

\smallskip

Hassett \cite{hassett} has studied the period map of cubic fourfolds in relation to the question whether a smooth cubic fourfold is rational or not. A key aspect of his work is the study of the {\it special cubic fourfolds}, i.e. the cubic fourfolds that acquire additional algebraic classes in $H^4(X)$ (N.B. since the Hodge conjecture holds for cubic fourfolds, this is equivalent to $\sH^{2,2}(X)\cap \sH^4(X,\bZ)$ having rank at least $2$). In particular, one considers the following notion:

\begin{definition}\label{defdethyp}
Fix a lattice $\Lambda=\langle1\rangle^{\oplus 21}\oplus \langle-1\rangle^{\oplus 2}$ and an element $h\in \Lambda$ of square $3$ with even orthogonal complement. Let $\Lambda_0=\langle h\rangle^{\perp}_{\Lambda}\cong E_8^{\oplus 2}\oplus U^{\oplus 2}\oplus A_2 $, and $\calD$ and $\Gamma$ as above. For a rank $2$ lattice $M\hookrightarrow \Lambda$ primitively embedded in $\Lambda$ with  $h\in M$ define the hyperplane $\calD_M$ by
$$\calD_M=\{\omega\in \calD\mid \omega\perp M\}.$$
We say $\calD_M$ is a hyperplane of determinant $d=\det(M)$. Note that indeed $\calD_M$ is the restriction of a hyperplane $H_M\subset \bP(\Lambda_0\otimes_\bZ \bC)$ to $\calD$. The  hyperplanes of a given determinant form an arithmetic arrangement (w.r.t. $\Gamma$) of hyperplanes in the domain $\calD$.  We denote the arrangements of hyperplanes of determinants $2$ and $6$ by $\calH_{\infty}$ and $\calH_{\Delta}$ respectively.
\end{definition} 

As noted below, the hyperplanes $\calD_M$ of a given determinant $d$ form a single $\Gamma$-orbit (cf.  \cite[Prop. 3.2.2, Prop. 3.2.4]{hassett}). In particular, $\calH_\infty/\Gamma$ and $\calH_\Delta/\Gamma$ are irreducible hypersurfaces in $\calD/\Gamma$.
\begin{lemma} With notations as in Def. \ref{defdethyp}, the following holds:
\begin{itemize}
\item[i)] $\calD_M$ is non-empty if and only if $\det(M)\equiv 0,2 \mod 6$;
\item[ii)] if $M$ and $M'$ have the same determinant, then $\calD_M$ and $\calD_{M'}$ are conjugated by $\Gamma$. 
\end{itemize}
\end{lemma}

The relevance for us of the arrangements $\calH_\infty$ and $\calH_\Delta$ comes from the observation of Hassett \cite[\S4.4]{hassett} (see also \cite[pg. 596, Prop. 1]{voisin}) that  the period map for smooth cubic fourfolds misses these two hyperplane arrangements. Specifically, the following holds:
\begin{proposition}\label{badhyper} 
Let $X$ be a smooth cubic fourfold. Then, $\sH^4_0(X,\bZ)\cap \sH^{2,2}(X)$  does not contain a primitive class $\delta$ such that either
\begin{itemize}
\item[i)] $\delta^2=2$,
\item[ii)] or $\delta^2=6$ and $\delta.x\equiv 0 \mod 3$ for all $x\in \sH^4_0(X,\bZ).$
\end{itemize}
The hyperplanes orthogonal to $\delta$ are the hyperplanes of determinant $6$ and $2$ respectively. Consequently, the image of the period map satisfies: $\textrm{Im}(\calP_0)\subset (\calD\setminus (\calH_{\infty}\cup \calH_{\Delta}))/\Gamma$.  
\end{proposition}

We recall that a primitive element $\delta$ in a lattice $L$ is called a {\it generalized root} iff $\delta^2\mid 2 \delta. x$ for all $x\in L$. For such an element $\delta$, the reflection $s_\delta$, given by $s_\delta(x)=x-2\frac{(x,\delta)}{(\delta,\delta)} x$, is an element of the orthogonal group $O(L)$. In our situation, both cases occurring in the previous proposition correspond to generalized roots in $\Lambda_0$, of norm $2$ and $6$ respectively. Consequently, we define: 

\begin{definition}
We call an element $\delta\in \Lambda_0$ satisfying $\delta^2=2$  a {\it root} of $\Lambda_0$, and an element $\delta\in \Lambda_0$ with $\delta^2=6$  and $\delta.\Lambda_0\equiv 0 \mod 3$  a {\it long root} of $\Lambda_0$. 
\end{definition}

Since a root $\delta$ corresponds geometrically to a vanishing cycle, the subgroup of $O(\Lambda_0)$ generated by reflections $s_\delta$ in the roots $\delta$ of $\Lambda_0$ coincides with the monodromy group $\Gamma\cong \mathrm{O^*}(\Lambda_0)$. A reflection $s_{\delta}$ in a long root $\delta$ acts non-trivially on the discriminant group $A_{\Lambda_0}\cong \bZ/3$ by  switching the generators. Thus, $s_{\delta}\in \mathrm{O^+}(\Lambda_0)$, but $s_\delta\not \in \mathrm{O^*}(\Lambda_0)\cong \Gamma$ (N.B. $\mathrm{O^+}(\Lambda_0)\cong \{\pm \mathrm{id}\}\times  \mathrm{O^*}(\Lambda_0)\cong \bZ/2\times  \mathrm{O^*}(\Lambda_0)$). It follows that the hyperplanes $\calH_{\Delta}$ are precisely the reflection hyperplanes for  $\Gamma$ (the hyperplanes in  $\Lambda_0\otimes_\bZ \bC$ pointwise fixed by some element of $\Gamma$). Similarly, the hyperplanes in $\calH_{\Delta}\cup \calH_\infty$ are the reflection hyperplanes for   $\mathrm{O^+}(\Lambda_0)$. Note however that the period domain is a subset of $\bP(\Lambda_0\otimes_\bZ \bC)$. In particular, we get that the hyperplanes from $\calH_{\infty}$ are anti-invariant hyperplanes in $\Lambda_0\otimes_\bZ \bC$, but become (pointwise) invariant in $\calD$ with respect to some $s_{\delta}$, where $\delta$ is a long root (N.B. by abuse of notation  we use $\calH_{\infty}$ to refer to hyperplanes both in $\Lambda_0\otimes_\bZ \bC$ and $\calD\subset \bP(\Lambda_0\otimes_\bZ \bC)$). 

\begin{lemma}\label{lemmaloop}
Let $\delta\in \Lambda_0$ be a long root. Then $(-s_\delta)\in O^*(\Lambda_0)$ and, via the isomorphism $O_h^+(\Lambda)\to O^*(\Lambda_0)$,  $(-s_\delta)$ acts as follows on $\Lambda$:
\begin{itemize}
\item[i)] it acts trivially on a primitive sublattice $M\cong \left(\begin{array}{cc} 3&2\\
2&2\end{array}\right)$ of $\Lambda$, with $M$ containing $h$ and $\delta$;
\item[ii)] it acts as multiplication by $-1$ on $M^\perp_\Lambda\cong A_1\oplus E_8^{\oplus 2}\oplus U^{\oplus 2}$.
\end{itemize}
\end{lemma}
\begin{proof}
Since $\delta.\Lambda_0\equiv 0 \mod 3$, it follows that  $h+\delta$  is divisible by $3$ in $\Lambda$ (possibly after changing $\delta$ to $-\delta$). The sublattice $M$ of $\Lambda$ generated by $h$ and $\frac{1}{3}(h+\delta)$ obviously satisfies the conditions of the lemma. For reasons that will be apparent later, we chose generators $h$ and $h-\frac{1}{3}(h+\delta)$ for $M$.
\end{proof}

\begin{remark}\label{remloop}
The element $\gamma=-s_\delta\in \Gamma$ (where $\delta$ is a long root) is  a Picard--Lefschetz type transformation associated to degenerations to the secant variety $X_0$ of the Veronese surface. Specifically, $\gamma$ is the monodromy transformation associated to a loop around the exceptional divisor of the Kirwan blow-up of the point $\omega$ (corresponding to $X_0$) in $\overline{\calM}$ (see section \ref{sectresolution}). The invariant cohomology under $\gamma$ (i.e. the lattice $M$) is the cohomology that comes from $X_0$. The anti-invariant cohomology is the vanishing cohomology, i.e.  the cohomology that it is determined by the direction of the degeneration to $X_0$. Furthermore, the Hodge structure supported on the vanishing cohomology is essentially the Hodge structure of a degree two $K3$ surface  (see section \ref{sectmodular}, esp. \ref{thmmonodromy2}). Indeed,  $M^\perp_\Lambda$ is (up to sign) the lattice associated to the primitive cohomology of a degree two $K3$ surface.
\end{remark}
%Further properties of the arrangements $\calH_{\Delta}$ and $\calH_\infty$ are discussed in section \ref{sectapplication}. 

%%%%%%%%%%%%%%%%%%%%%%%%%%%%%%%%%%%%%%%%%%%%%%%%%%%%%%%%%%%%%%
\subsection{Degenerations of Hodge Structures for cubic fourfolds}\label{sectdeg}
One of the main aspects of our work is the study of $1$-parameter families of cubic fourfolds from the point of view of variations of Hodge structures. In this section we collect a series of general results of this theory.
\begin{notation}
{\it A $1$-parameter degeneration} $f:\calX\to \Delta$ is a proper analytic map with generic fiber a smooth projective variety (where $\calX$ is an analytic variety and $\Delta$ is the unit disk). We denote by $X_0$ the special fiber, by $X_t$ the generic fiber, and by $X_\infty:=\calX^* \times_{\Delta^*}\calh$ the canonical fiber, where $\calX^*$ is the restriction of $\calX$ to the punctured disk $\Delta^*$ and  $\calh\to \Delta^*$ is the universal cover of $\Delta^*$ ($\calh$ is the upper half plane). We also say that $f$ is {\it a $1$-parameter smoothing of $X_0$}. 
\end{notation}

\subsubsection{Limit Mixed Hodge Structures} To a $1$-parameter degeneration there is associated a canonical limit mixed Hodge structure on $\sH^n_{\lim} :=\sH^n(X_\infty,\bQ)$ by Schmid and Steenbrink. To fix the notation and terminology we recall a few general facts (for a survey see \cite{clemensschmid}). The monodromy $T$ of the family over the punctured disk is quasi-unipotent. After a ramified base change of type $t\to t^r$ one can further assume that $T$ is actually unipotent. We will assume this in what follows. We fix $n$ to be the dimension of the fiber.  

Let $N=\log T$ be the logarithm of the monodromy. It is a nilpotent endomorphism acting on the cohomology $\sH^n_{\lim}$ with index of nilpotence $\nu$ (i.e. $N^\nu=0$, but $N^{\nu-1}\neq 0$). On the vector space $\sH^n_{\lim}$ there are two filtrations: a decreasing limit Hodge filtration $F^p$ and an increasing weight filtration $W_k$ induced by  $N$ (see \cite[pg. 106-107]{clemensschmid}). Together the two filtrations define a mixed Hodge structure on $\sH^n_{\lim}$. We are interested in determining the index of nilpotency $\nu$ of a $1$-parameter degeneration. For this, we note the relation:
$$\nu=\max \{k\mid\Gr^W_{n-k}\sH^n_{\lim}\neq 0\}+1.$$
We then note that the possibilities for the  Hodge numbers $h^{p,q}:=\dim \Gr_F^p\Gr_{p+q}^W(\sH^n_{\lim})$ are quite restricted. Namely, $\sH^n_{\lim}$ is isomorphic as a vector space to $\sH^n(X_t)$ and the following holds:
\begin{itemize}
\item[i)] $\Gr^W_k(\sH^n_{\lim})$ carries a pure Hodge structure (induced by the filtration $F^p$) of weight $k$;
\item[ii)] $N^k:\Gr^W_{n+k}(\sH^n_{\lim})\to\Gr^W_{n-k}(\sH^n_{\lim})$ is an isomorphism of Hodge structures of type $(-k,-k)$;
\item[iii)] $\dim_\bC F^p \sH^n_{\lim}=\dim_\bC F^p \sH^n(X_t)$ (see \cite[Cor. 11.25]{peterssteenbrink}).
\end{itemize}

\subsubsection{The case of cubic fourfolds} We specialize the above discussion to the case of degeneration of cubic fourfolds. First, since the Hodge structure on the middle cohomology of a smooth cubic fourfold is of level $2$ (i.e. $\sH^{p,q}=0$ if $|p-q|>2$), it follows that the index of nilpotence of $N$ is at most $3$.  Similarly to the case of K3 surfaces, we then define $3$ types of degenerations:
\begin{definition}\label{typedeg}
Let $f:\calX\to \Delta$ be a $1$-parameter degeneration of cubic fourfolds. Assume that the monodromy $T$ is unipotent.  We say that $f$ is a {\it Type I (II, or III) degeneration} if the index of nilpotence $\nu$ of $N=\log T$ is $1$ ($2$, or $3$ respectively).
\end{definition}

A simple observation is that, for cubic fourfolds, to determine the type of the degenerations it suffices to know the non-vanishing of a graded piece of $\sH^4_{\lim}$.

\begin{lemma}\label{lemmamonodromy}
Let $f:\calX\to \Delta$ be a $1$-parameter degeneration of cubic fourfolds. Then $f$ is a Type II degeneration iff $\Gr_3^W\sH^4_{\lim}\neq 0$. In this case $\Gr_3^W\sH^4_{\lim}$ is a Tate twist of the Hodge structure of an elliptic curve $E$, i.e. $\Gr_3^W\sH^4_{\lim}\cong \sH^1(E)(-1)$. Similarly,  $f$ is a Type III degeneration iff $\Gr_2^W\sH^4_{\lim}\neq 0$, in which case $\Gr_2^W\sH^4_{\lim}$ is a trivial $1$-dimensional Hodge structure of weight $2$.
 \end{lemma}
 \begin{proof}
 Since $X_t$ is  a smooth cubic fourfold, we have $\dim_\bC F^4\sH^4(X_t)=0$ and $\dim_\bC F^3\sH^4(X_t)=1$. Thus, the only possibly non-zero Hodge numbers in weight at most  $3$ are $h^{2,1}$, $h^{1,2}$, and $h^{1,1}$. Additionally, they satisfy (by item iii above) $h^{3,1}+h^{2,1}+h^{1,1}=1$. The claim follows.
 \end{proof}

We note that the type of a degeneration is closely related to the Baily--Borel compactification $(\calD/\Gamma)^*$ as follows. A $1$-parameter degeneration of cubic fourfolds induces a period map  $g:\Delta^*\to \calD/\Gamma$, which always extends to an analytic map $\Delta\to  (\calD/\Gamma)^*$. For Type I degenerations the limit point $\lim_{z\to 0} g(z)$ belongs to the interior $\calD/\Gamma$. In contrast, the limit point of a Type II (Type III) degeneration belongs to a Type II (Type III respectively) boundary component (see \S\ref{sectbb}). 

\subsubsection{The specialization morphism} The formalism of vanishing cycles (\cite[Exp. XIII]{sga72}) relates the limit mixed Hodge structure of a smoothing  to the mixed Hodge structure of the central fiber via the specialization morphism:
$$\textrm{sp}_n:\sH^n(X_0)\to \sH^n_{\lim}.$$
We recall the basic construction as needed in our situation. 

\smallskip

First, from the specialization diagram:
\begin{equation}\label{specializationdiag}
\xymatrix@R=.25cm{
{X_0}\ar@{->}[r]^{i}&{\calX} \ar@{->}[dd]_{f} \ar@{<-}[rr]^{\pi} &&{X_{\infty}}\ar@{->}[dd] \\
&&\\
&{\Delta}\ar@{<-}[r]&{\Delta^*}\ar@{<-}[r]&{\calh}
}
\end{equation}
we define {\it the functor of nearby cycles} $\psi_f:\mathrm{D}_c^b(\calX)\to \mathrm{D}_c^b(X_0)$ (where $D_c^b(.)$ denotes the bounded derived category of constructible complexes) by 
$$\psi_f \calF^{\bullet}:=i^*R\pi_*\pi^*\calF^{\bullet}$$ 
(see \cite[Def. 4.2.1]{dimca}). There exists a natural comparison map $i^*\calF^{\bullet}\xrightarrow{c}\psi_f \calF^{\bullet}$, and the {\it specialization morphism} $\textrm{sp}_n$ is defined to be the cohomology map associated to $c$. The {\it functor of vanishing cycles} $\phi_f$ is the cone over the morphism $c$. By definition there exists a distinguished triangle:
 \begin{equation}\label{triangle}
\xymatrix@R=.25cm{
{i^*\calF^{\bullet}}\ar@{->}[rr]^{c}&&{\psi_f \calF^{\bullet}} \ar@{->}[ddl]^{\mathrm{can}} \\
&&\\
&{\phi_f \calF^{\bullet}}\ar@{->}[uul]^{[1]}&
}
\end{equation}
in derived category  $\mathrm{D}_c^b(X_0)$ (see \cite[\S4.2]{dimca}). 

\smallskip

We are interested in the situation when $\calF^{\bullet}$ is the constant sheaf $\underline{\bC}_\calX$. Taking the hypercohomology associated to (\ref{triangle}), we obtain a long exact sequence  relating the cohomology of the central fiber with the cohomology of the canonical fiber $X_{\infty}$:
\begin{equation}\label{specialization}
\dots\to \sH^n(X_0)\xrightarrow{\textrm{sp}_n} \sH^n_{\lim}\to \bH^n(\phi_f \underline{\bC}_\calX)\to\sH^{n+1}(X_0)\to \dots
\end{equation} 
 Furthermore, {\it the vanishing cohomology} $\bH^n(\phi_f \underline{\bC}_\calX)$ can be endowed with a natural mixed Hodge structure making (\ref{specialization}) an exact sequence of mixed Hodge structures.  Since a morphism of mixed Hodge structures is strict with respect to both the weight and Hodge filtration, $\mathrm{sp}_n$ maps $\sH^{p,q}(X_0)$ to $\sH^{p,q}_{\lim}$ (where $\sH^{p,q}=\Gr_F^p\Gr^W_{p+q}\sH$ for a mixed Hodge structure $\sH$). In particular, the statement {\it $\textrm{sp}_n$ is an isomorphism on the $(p,q)$ components} is well-defined.  

\smallskip

In our situation, $X_t$ is either a cubic fourfold or a $K3$ surface, and $n$ represents the dimension of $X_t$. In particular, there is no odd cohomology. Thus, we have:
\begin{equation}\label{assumption}
\sH^{n-1}(X_t)=\sH^{n+1}(X_t)=0,
\end{equation} 
It follows that (\ref{specialization}) reduces to a five-term exact sequence:
\begin{equation}\label{exactsequence}
0\to\bH^{n-1}(\phi_f \underline{\bC}_\calX)\to \sH^n(X_0)\xrightarrow{\textrm{sp}_n} \sH^n_{\lim}\to \bH^n(\phi_f \underline{\bC}_\calX)\to\sH^{n+1}(X_0)\to 0.
\end{equation} 

To compute the vanishing cohomology we note the following spectral sequence: 
\begin{equation}\label{leray}
E^{p,q}_2=\sH^p(X_0,\calH^q(\phi_f \underline{\bC}_\calX))\Longrightarrow \bH^{p+q}(\phi_f \underline{\bC}_\calX)
\end{equation}
Furthermore,  the stalk of the cohomology sheaf $\calH^q$ is the reduced cohomology of the Milnor fiber:
\begin{equation}\label{stalk}
\calH^q(\phi_f \underline{\bC}_\calX)_x=\widetilde{\sH}^{q}(F_x;\bC),
\end{equation}
where as usually $F_x$ is the intersection of the generic nearby fiber with  a small open ball centered at $x$.
% (see \cite[pg. 106]{dimca}). 
In particular, $\calH^q(\phi_f \underline{\bC}_\calX)$ are supported on the singular locus of $X_0$ (\cite[Prop. 4.2.8]{dimca}). The total space  $\calX$  of the degeneration might be singular, but in our situation (degeneration of hypersurfaces) it has at worst local complete intersection singularities. Therefore, the range for which there exists non-vanishing cohomology for the Milnor fiber is the same as in the smooth case (\cite[Prop. 6.1.2]{dimca}). 

\smallskip

In particular, in the case that the special fiber $X_0$ has at worst isolated singularities, one obtains that the specialization morphism $\textrm{sp}_n$ is injective (\cite[Prop. 2.7]{dolgachevinsignificant},  \cite[6.2.2, 6.2.4]{dimca}) and an exact sequence of mixed Hodge structures:
\begin{equation}\label{spectralisolated}
0\to \sH^n(X_0)\xrightarrow{\textrm{sp}_n} \sH^n_{\lim}\to \oplus_{x_i\in \Sing(X_0)} \sH^n(X_i)\to \sH^{n+1}(X_0)\to 0
\end{equation} 
where $\sH^n(X_i)$ is the cohomology of the Milnor fiber at $x_i$ endowed with the mixed Hodge structure of Steenbrink \cite{steenbrinkmhs}. Similarly, if the central fiber has $1$-dimensional locus the first part of the exact sequence (\ref{exactsequence}) reads:
\begin{equation}\label{spectralnonisolated}
0\to \sH^0(X_0,\calH^{n-1}(\phi_f \underline{\bC}_\calX))\to \sH^n(X_0)\xrightarrow{\textrm{sp}_n} \sH^n_{\lim}\to \dots
\end{equation}

%%%%%%%%%%%%%%%%%%%%%%%%%%%%%%%%%%%%%%%%%%%%%%%%%%%%%%%%%%%%%%
\section{The monodromy around semistable cubic fourfolds}\label{sectmonodromy}
The main result of this section is the control of the monodromy for degenerations of cubic fourfolds which are not of Type IV.
\begin{theorem}\label{thmmonodromy}
Let $f:\calX\to \Delta$ be a $1$-parameter smoothing of a semi-stable cubic fourfold $X_0$ with closed orbit.  The following holds:
\begin{itemize}
\item[i)] if $X_0$ has Type I, then $f$ is a Type I degeneration;
\item[ii)] if $X_0$ has Type II, then $f$ is a Type II degeneration;
\item[iii)] if $X_0$ has Type III, then $f$ is a Type III degeneration.
\end{itemize}
(see definitions \ref{deftypecubic} and \ref{typedeg}).
\end{theorem}
\begin{proof}
By lemma \ref{lemmamonodromy}, the monodromy of the family is determined by the non-vanishing of  $\sH_{\lim}^{p,1}$ for $p=3$, $2$, or $1$. By proposition \ref{reductioncentralfiber}, the non-vanishing of  $\sH^{p,1}(X_0)$ implies the non-vanishing of $\sH_{\lim}^{p,1}$. The claim now follows from the computation, done in proposition \ref{computetype23}, of the mixed Hodge structure of the central fiber.
\end{proof}
 
Due to the finiteness of the monodromy for Type I fourfolds a stronger result holds. Namely,  the period map extends over the locus $\calM$ of such fourfolds.

\begin{proposition}\label{extension1}
The period map for a cubic fourfolds $\calP_0:\calM_0\to \calD/\Gamma$ extends to a morphism $\calP:\calM\to \calD/\Gamma$ over the simple singularity locus $\calM\subset \overline{\calM}$. The image of $\calM\setminus \calM_0$ under the extended period map $\calP$ is contained in $\calH_{\Delta}/\Gamma$.
\end{proposition} 
\begin{proof}
Let $o\in \calM\setminus \calM_0$ correspond to a cubic fourfold $X_0$ with simple isolated singularities. The statement is analytically local at $o$, and stable by finite base changes. Since $\calM$ is a geometric quotient, after shrinking and  a possible finite cover, we can assume that a neighborhood of $o$ in $\calM$ is a $20$-dimensional ball $S$. We can further assume that there exists a family of cubic fourfolds $\calX\to S$ with at worst simple isolated singularities and fiber $X_0$ over $o$. Let $o\in \Sigma$ be   discriminant hypersurface. Over $S\setminus \Sigma$ the family $\calX$ gives a variation of Hodge structures defining the period map $\calP_0$. By the removable singularity theorem \cite[pg. 41]{griffithscurvature}, the extension statement is equivalent to  the monodromy representation $\pi_1(S\setminus \Sigma,t)\to \mathrm{Aut}(\sH^4(X_t,\bZ))$ (for $t\in S\setminus \Sigma$) having finite image $\Gamma_0$ (N.B. $\Gamma_0$ is the local monodromy group around $X_0$).

The fourfold $X_0$ is GIT stable, thus it has finite stabilizer. For cubic fourfolds we have:
$$(*)\ \ \  \ n(d-2)-2=d,$$ 
where $n-1=4$ is the dimension and $d=3$ the degree. It follows that the conditions (i.e. finite stabilizer and $(*)$) from du Plessis--Wall  \cite[Cor. 1.6]{versality} and \cite{duplessiswall0} are satisfied. Thus,  the family $\calX$ gives a  simultaneous versal deformation of the singularities of $X_0$. It follows that the local monodromy group $\Gamma_0$ is the product of the monodromy groups associated to the singularities of $X_0$. Since the singularities of $X_0$ are of type A-D-E and the dimension is even,  these monodromy groups are (finite) Weyl groups of type A-D-E. Thus, as needed, $\Gamma_0$ is finite. Furthermore, the period point corresponding to $X_0$ is left invariant by the reflections in the vanishing cycles. Thus,  it belongs to $\calH_{\Delta}/\Gamma$. 
\end{proof}

%%%%%%%%%%%%%%%%%%%%%%%%%%%%%%%%%%%%%%%%%%%%%%%%%%%%%%%%%%%%%%
\subsection{Reduction to the central fiber}\label{reductioncentral}
\begin{proposition}\label{reductioncentralfiber}
Let $X_0$ be a semistable cubic fourfold with closed orbit of Type I--III. Let $\calX\to \Delta$ be any $1$-parameter smoothing of $X_0$,  and consider the associated specialization morphism $\textrm{sp}_4:\sH^4(X_0)\to \sH^4_{\lim}$. If $X_0$ has isolated (non-isolated) singularities then $\textrm{sp}_4$  induces an isomorphism (resp. injection) on the $(p,q)$ components of corresponding mixed Hodge structures for all $p$ and $q$ with $p+q\le 4$ and $(p,q)\neq (2,2)$. 
\end{proposition}
\begin{proof} We divide the proof in two cases: either $X_0$ has isolated singularities, or not. 

\smallskip

\noindent{\bf Case 1 (Isolated singularities):} Assume that $X_0$ has only isolated singularities which are suspensions of the types listed in definition \ref{shahlist}. From the exact sequence (\ref{spectralisolated}):
\begin{equation*}
0\to \sH^4(X_0)\xrightarrow{\textrm{sp}_4} \sH^4_{\lim}\to \oplus_{x_i\in \Sing(X_0)} \sH^4(X_i) \dots
\end{equation*} 
and the strictness of the morphisms of mixed Hodge structures, we see that it is enough to prove the following claim:

\smallskip

\noindent$(*)$ {\it If $\Gr_F^p\Gr_{p+q}^W \sH^4(X_i)\neq 0$ then $(p,q)\in \{(2,2), (3,2), (2,3), (3,3)\}$, where $\sH^4(X_i)$ is the vanishing coholomogy of a smoothing of a simple, simple elliptic, or cusp singularity.} 

\smallskip

We note the following facts about the mixed Hodge structure on the vanishing cohomology:
\begin{itemize}
\item[i)] ({\it the local nature}) the mixed Hodge structure on $\sH^4(X_i)$ depends only on the germ of the smoothing $(\calX,x_i)$ (cf. \cite[pg. 560]{steenbrinksemicontinuity});
\item[ii)] ({\it semicontinuity}) $\dim \Gr_F^p \sH^4(X_i)$ is independent of the smoothing (cf. \cite[Cor. 2.6]{steenbrinksemicontinuity}).
\end{itemize}
From the semicontinuity property, it follows immediately that if $(*)$ holds for a smoothing, then it holds for any smoothing. We therefore check $(*)$ for the standard Milnor fibrations of the singularities $A_n$, $D_m$, $E_r$, $\widetilde{E_r}$, and $T_{p,q,r}$ respectively. For these singularities the computation of the mixed Hodge structure on the Milnor fiber is well known (see \cite[II.8]{kulikovmhs}). In fact, $(*)$ is equivalent to the statement that the spectrum of those singularities (in dimension $4$) is included in the interval $[1,2]$ (see also remark \ref{remspectrum}). This settles the case of isolated singularities. 

\medskip

\noindent{\bf Case 2 (Non-isolated singularities):} Let $C$ be the $1$-dimensional singular locus of $X_0$, and $C_i$ its  irreducible components. According to proposition \ref{propsingtype},  the singularities along $C$ are of type $A_{\infty}$ at all but a finite number of points. At these special points the singularities are of type $D_{\infty}$ or degenerate cusps (non-isolated singularities of type $t_4$--$t_6$, see \ref{shahlist}). We denote by $\mathring{C}$ (and $\mathring{C_i}$) the non-special locus. We note two simple facts:
\begin{itemize}
\item[i)] $\mathring{C}=\cup_i \mathring{C_i}$ is disconnected (i.e. any point of intersection is special);
\item[ii)] $\pi_1( \mathring{C_i})$ is non-trivial (in fact, either $C_i$ is elliptic and $ \mathring{C_i}=C_i$, or $C_i$  is rational and there are either $2$ or $4$ special points depending on the type of the degeneration).
\end{itemize}

We recall, the exact sequence (\ref{spectralnonisolated}):
\begin{equation}
0\to \sH^0(X_0,\calH^{3}(\phi_f \underline{\bC}_\calX))\to \sH^4(X_0)\xrightarrow{\textrm{sp}_4} \sH^4_{\lim}\to \dots
\end{equation}  
and that 
$$\calH^{3}(\phi_f \underline{\bC}_\calX)_x\cong \sH^3(F_x,\bC)$$
where $F_x$ is the Milnor fiber at $x$. In particular, the sheaf of vanishing cycles is supported on the $1$-dimensional singular locus. The injectivity statement is equivalent to saying that  the mixed Hodge structure on the vanishing cohomology satisfies  $\Gr_F^p\sH^0(X_0,\calH^{3}(\phi_f \underline{\bC}_\calX))=0$ for $p<2$.

\smallskip

The analysis of the vanishing cycles in the case of $1$-dimensional singular locus is relatively well understood (see \cite{siersma}, \cite[II.8.10]{kulikovmhs}). Specifically, we stratify $C=\Sigma_1\cup \Sigma_0$ such that the vanishing cycles form a local system over $\Sigma_1$. In our situation, $\Sigma_1=\cup_i  \mathring{C_i}$, and $\Sigma_0$ are the special points. Since, the transversal singularity is $A_1$, we obtain a $1$-dimensional local system over $\Sigma_1$.  Since the vertical monodromy, given by the natural action of $\pi_1(\mathring{C_i})$  (see \cite[Ch. 3]{siersma}), is non-trivial on each component, there are no non-zero sections of the local system over $\Sigma_1$. Thus,  the sections of $\calH^3(\phi_f \underline{\bC}_\calX))$ are supported on the special points. In the type II case  we only have special points of type $D_{\infty}$, but then $\sH^3(F_x)=0$ (\cite[pg. 5]{siersma}). Thus, $\sH^0(X_0,\calH^{3}(\phi_f \underline{\bC}_\calX))=0$ for the Type II case. For the Type III case, we note that the situation is local around the special points. Via the generalized Thom-Sebastiani theorem  \cite{saitosc,loeser} the computation of the vanishing cohomology is reduced to  the surface case, but since the singularities  are restricted to the list given in \ref{shahlist} we obtain $\Gr_F^p\sH^0(X_0,\calH^{3}(\phi_f \underline{\bC}_\calX))=0$ for $p<2$ (see \cite[Thm. 2]{shahinsignificant}). 
\end{proof}

\begin{remark}\label{remspectrum}
 We note that the analogue of $(*)$ in the surface case is equivalent to the condition that spectrum is included in $[0,1]$ (N.B. the spectrum is translated by $\frac{1}{2}$ for each suspension). This in turn is equivalent to the singularities being log canonical. In conclusion, the isolated singularities for which proposition \ref{reductioncentralfiber}  holds are precisely the simple, simple-elliptic, and cusp singularities (see also remark \ref{reminsignficant},  \cite[Thm. 1.2]{shahinsignificant}, and \cite[Thm. 4.13]{dolgachevinsignificant}).
\end{remark}

%%%%%%%%%%%%%%%%%%%%%%%%%%%%%%%%%%%%%%%%%%%%%%%%%%%%%%%%%%%%%%
\subsection{The Mixed Hodge Structure of  Type I--III fourfolds}\label{sectmhs} In this section we discuss the mixed Hodge structure on the middle cohomology of a singular cubic fourfold. The basic observation is that a singular cubic fourfold $X_0$ is rational: $X_0$ is birational to $\bP^4$ via the projection from any singular point $p$. This allows us to reduce the computation of the mixed Hodge structure on $\sH^4(X_0)$  to a computation for degenerations of K3 surfaces. We obtain the following result, concluding the proof of Theorem \ref{thmmonodromy}.
  
\begin{proposition}\label{computetype23}
Let $X_0$ be a semistable cubic fourfold with closed orbit. Then the following holds:
\begin{itemize}
\item[i)]  if $X_0$ is of Type I, then $\sH^4(X_0)$ carries a pure  Hodge structure of weight $4$ with $\sH^{3,1}(X_0)\neq 0$;
\item[ii)] if $X_0$ is of Type II, then $\Gr_3^W\sH^4(X_0)\neq 0$;
\item[iii)] if $X_0$ is of Type III, then $\Gr_2^W\sH^4(X_0)\neq 0$.
\end{itemize}
\end{proposition}
\begin{proof}
If $X_0$ is smooth there is nothing to prove. Assume that $X_0$ is singular. Let $p\in \Sing(X_0)$ (suitably chosen) and $\pi_p:X_0\dasharrow \bP^4$ the projection map. In \S\ref{sectrational} (esp. Cor. \ref{reduction3}), we relate the Hodge structure of $X_0$ to the Hodge structure of the surface $S_p$, the base locus of $\pi_p^{-1}$. The proposition then follows from the  analysis of the mixed Hodge structure on $\sH^2(S_p)$. But since the surface $S_p$ is a degeneration of K3 surfaces with insignificant limit singularities, the computation is standard. The case by case analysis is done in \ref{lemmacase1},  \ref{lemmacase2}, and \ref{lemmacase3} for Type I-III respectively.
\end{proof}

\subsubsection{Singular cubic fourfolds and degree $6$ $K3$ surfaces}\label{sectrational} Let $X_0$ be a singular cubic fourfold (irreducible and reduced). We choose a singular point $p\in X_0$ and  assume additionally:
\begin{center}
{\it (*) $\corank_p(X_0)\le 3$ and no line contained in $\Sing(X_0)$ passes through $p$.} 
\end{center}
The linear projection $\pi_p$ with center $p$ gives a birational isomorphism between $X_0$ and $\bP^4$. The birational map $\pi_p$ can be resolved by blowing-up the point $p$. The result is the following diagram:
\begin{equation}\label{blowupdiagram}
\xymatrix{
&Q_p \ar@{->}[ld]\ar@{^{(}->}[r]&\widetilde{X_0}\ar@{<-^{)}}[r] \ar@{->}[ld]_{f} \ar@{->}[rd]^{g}&E_p\ar@{->}[rd] \\\
p\ar@{^{(}->}[r]&X_0 \ar@{-->}[rr]^{\pi_p}&&\bP^4 \ar@{<-^{)}}[r] &S_p 
}
\end{equation}
where  $Q_p$ is the projectivized tangent cone at $p$,  $S_p\subset \bP^4$ is the surface parametrizing the lines of $X_0$ through $p$, and $E_p=g^{-1}(S_p)$ is the exceptional divisor of $g$. The following facts about the surface $S_p$ are  well known (see esp. O'Grady \cite[\S5.4]{ogrady}).

\begin{proposition}\label{propblowup}
With notations as above, assume that $p\in X_0$ is a singular point satisfying (*). Then the following holds:
\begin{itemize}
\item[i)] $S_p$ is the  complete intersection of  a quadric and cubic in $\bP^4$; 
\item[ii)] $S_p$ is reduced (but possibly reducible);
\item[iii)] $S_p$ has only hypersurface singularities;
\item[iv)] the singularities of $S_p$ are in one-to-one correspondence,  including the type,  with the singularities of $\widetilde{X_0}$.
\end{itemize}
Furthermore, the morphism $g$ of (\ref{blowupdiagram}) is the blow-up of $\bP^4$ along $S_p$. In particular,  $E_p$ is a $\bP^1$-bundle over $S_p$.
\end{proposition}
\begin{proof}
The first statement is well-known (e.g. \cite[Remark 5.11]{ogrady}). Namely, let $X_0$ be given by:
$$X_0:\ \  (x_0Q(x_1,\dots,x_5)+F(x_1,\dots,x_5)=0),$$
with $Q$ and $F$ non-vanishing homogeneous polynomials of degree $2$ and $3$ respectively. Then, the surface $S_p$ is the complete intersection of $Q$ and $F$. The two assumptions from (*) imply that $Q$ is reduced and that $Q$ and $F$ are not simultaneously singular. In particular, since at any point of $S_p$ we can choose either $Q$ or $F$ as a local coordinate, we obtain iii). The relation between the singularities of $S_p$ and those of $X_0$ and $\widetilde{X_0}$ was analyzed in O'Grady \cite[Prop. 5.15]{ogrady} and Wall \cite[\S I.2]{wallsextic}. In particular, it is not hard to see that, under the assumption (*), the singularities of $\widetilde{X_0}$ are double suspensions of the singularities of $S_p$  (see \cite[pg. 7]{wallsextic}), i.e. the singularities of $\widetilde{X_0}$ and $S_p$ have the same analytic type.

The  statement about $g$ is  \cite[Prop. 5.14]{ogrady}. Finally, the exceptional divisor $E_p$ is the projectivization of the normal bundle $\calN_{S_P/\bP^4}$. Since $S_p$ is a complete intersection, $\calN_{S_P/\bP^4}$  is a rank $2$ vector bundle over $S_p$ (see \cite[App. B \S6--7]{fultonintersection}) and the claim follows.\end{proof}

Since both $f$ and $g$ are explicit blow-ups, we are able to compute the Hodge structures of $X_0$ from that of $S_p$. Specifically, we recall that for a proper birational modification there exists a  Mayer-Vietoris sequence relating the cohomologies of the spaces involved: Given a diagram 
$$
\xymatrix@R=.25cm{
{E} \ar@{->}[dd] \ar@{->}[r]&{\widetilde{X}}\ar@{->}[dd] \\
&&\\
{D}\ar@{->}[r]&{X}
}
$$
with $X$ and $\widetilde{X}$ projective varieties, $f:\widetilde{X}\to X$ a projective birational morphism, $D$ the discriminant of $f$ and $E=\pi^{-1}(D)$, there exists a long exact sequence of mixed Hodge structures
\begin{equation}\label{mayervietoris}
\dots \to \sH^{n-1}(E)\to \sH^n(X)\to \sH^n(\widetilde{X})\oplus \sH^n(D)\to \sH^n(E)\to \sH^{n+1}(X)\to \dots
\end{equation}
(see \cite[Cor. 5.37]{peterssteenbrink}). Applying this exact sequence to the discriminant square associated to the morphism $f$ (see  (\ref{blowupdiagram})), we relate the cohomology of $X_0$ with that of the blow-up $\widetilde{X_0}$: 

\begin{lemma}\label{reduction1}
Let $X_0$ be a singular cubic fourfold (irreducible, reduced), and $p$ a singular point satisfying (*).  Then, the pullback $f^*:\sH^4(X_0)\to \sH^4(X_0)$ induces a $(p,q)$-isomorphism for all $(p,q)\neq (2,2)$.
\end{lemma}
\begin{proof}
The morphism $f$ is the blow-up of the singular point $p$. Thus, the relevant part of the sequence (\ref{mayervietoris}) reads:
$$ \dots \to \sH^{3}(Q_p)\to \sH^4(X_0)\to \sH^4(\widetilde{X_0})\to \sH^4(Q_p)\to \dots$$
Since for the reduced quadric $Q_p$ we have $\sH^3(Q_p)=0$ and $\sH^4(Q_p)$ is of pure weight $4$ and type $(2,2)$, the claim follows. 
\end{proof}

Similarly, we relate the cohomology of the $\widetilde{X_0}$ with that of the surface $S_p$:
\begin{lemma}\label{reduction2}
With notations and assumptions as above, 
$$\dim \Gr_F^1\Gr^W_{n+2}\sH^4(\widetilde{X_0})=\dim \Gr_F^0\Gr^W_n\sH^2(S_p)$$
for $n=0,1,2$.
\end{lemma}
\begin{proof}
We consider the exact sequence (\ref{mayervietoris}) for the birational morphism $g$ of (\ref{blowupdiagram}). We obtain:
$$ \dots \to \sH^{4}(\bP^4)\to \sH^4(\widetilde{X_0})\oplus \sH^4(S_p)\to \sH^4(E_p)\to \sH^5(\bP^4)=0$$
Both $\sH^4(\bP^4)\cong \bC$ and $\sH^4(S_p)$ carry a pure weight $4$ Hodge structure of type $(2,2)$ (for $S_p$ this follows from \cite[Thm. 6.33]{peterssteenbrink}). Thus, the restriction map $\sH^4(\widetilde{X_0})\to \sH^4(E_p)$ is a $( p,q)$-isomorphism for all $(p,q)\neq (2,2)$.  The claim now follows from the fact that $E_p$ is a projective  $\bP^1$-bundle over $S_p$ (cf. Prop. \ref{propblowup}).  \end{proof}

In conclusion, we obtain:
\begin{corollary}\label{reduction3}
Let $X_0$ be an irreducible, reduced cubic fourfolds, and $p\in X_0$ a singular point satisfying (*). Then $\Gr_F^0\Gr^W_{n+2}\sH^4(X_0)=0$ and 
\begin{equation}\label{eqreduction3}
\dim \Gr_F^1\Gr^W_{n+2}\sH^4(X_0)=\dim \Gr_F^0\Gr^W_n\sH^2(S_p)
\end{equation}
for $n=0,1,2$.
\end{corollary}

\subsubsection{The Type I case}\label{subcase1} 
The surface $S_p$ associated to a singular Type I fourfold $X_0$ has at worst du Val singularities. Thus, $S_p$ is a $K3$ surface. It follows that the Hodge structure on $\sH^4(X_0)$ is pure.

\begin{lemma}\label{lemmacase1}
Let $X_0$ be a cubic fourfold with at worst simple isolated singularities. Then $\sH^4(X_0)$ carries a pure Hodge structure of weight $4$. Furthermore, $\sH^{3,1}(X_0)\neq 0$. 
\end{lemma}
\begin{proof}
Assume that $X_0$ is singular, and choose any singular point $p\in X_0$. According to proposition \ref{propblowup}, the singularities of $S_p$ are in bijective correspondence to the singularities of $\widetilde{X_0}$ and of the same analytic type. The singularities of $\widetilde{X_0}$ are of two kinds: either they come from $X_0\setminus\{ p\}$, or they lie over $p$. In either case they are simple singularities. It follows that $S_p$ is a $(2,3)$-complete intersection in $\bP^4$ with at worst du Val singularities.  It is well known that the minimal desingularization of $S_p$ is a $K3$ surface and that the Hodge structure on $H^2(S_p)$ is pure. The lemma follows (cf. \ref{reduction3}). \end{proof}

\begin{remark}
An essentially equivalent statement is O'Grady \cite[Prop. 5.9, 5.28]{ogrady}. We remark that the proof of   \cite[Prop. 5.9]{ogrady} (see  \cite[\S5.4.4]{ogrady}) actually  works under the assumption that $S_p$ is a surface with at worst du Val singularities. This gives another proof to the lemma.
\end{remark}

\subsubsection{The Type II and III case}\label{subcase2} 
In the case of Type II and III fourfolds, we can reduce the computation of the mixed Hodge structure on $\sH^4(X_0)$ to a computation for degenerations of degree $6$ $K3$ surfaces. This reduction is made possible by the correspondence of singularities given by \ref{propblowup} and the fact that the singularities that occur are the insignificant singularities of Shah \cite{shahinsignificant}. 

\begin{lemma}\label{lemmacase2}
Let $X_0$ be a Type II cubic fourfold with closed orbit. Then we can chose a singular point $p\in X_0$ satisfying (*) such that the surface $S_p$ associated to the projection from $p$  is a Type II degeneration of K3 surfaces,  i.e. $S_p$ has insignificant limit singularities and $\Gr_1^W\sH^2(S_p)\neq 0$. Thus, $\Gr_3^W\sH^4(X_0)\neq 0$.
\end{lemma}
\begin{proof} 
All singularities of  a Type II fourfold $X_0$ have corank at most $3$. Also, in each of the cases $\alpha$--$\phi$  we can find a singular point $p$ not lying on a singular line (see table \ref{tableboundary}). Thus, we can choose a point $p$ satisfying (*). We can further assume that $p$ is either of type $\widetilde{E}_r$ or $A_{\infty}$. It is immediate to check that the non-simple singularities of the fourfold $\widetilde{X_0}$ (obtained by blowing-up $p$) are of type   $\widetilde{E_r}$, $A_{\infty}$, or $D_{\infty}$. By proposition \ref{propblowup}, the same holds for $S_p$. In particular, we conclude that $S_p$ has only insignificant limit singularities and is a degeneration of (degree $6$) $K3$ surfaces.  From Shah \cite[Thm. 2]{shahinsignificant} it follows that $\dim \Gr^0_F \sH^2(S_p)=1$. Since by construction $S_p$ has at least one non-du Val singularity, the Hodge structure on $S_p$ is not pure. Thus, $\Gr_k^W\sH^2(S_p)\neq 0$  for either $k=1$ (Type II) or $k=0$ (Type III). It is known that if all the singularities are of type $\widetilde{E_r}$, $A_{\infty}$, or $D_{\infty}$ the degeneration is of Type II (see \cite[Thm. 3.2]{shah}), concluding the proof of the lemma.

For concreteness, we illustrate the case $\delta$. In this case, the cubic fourfold $X_0$ has $3$ singularities of type $\widetilde{E_6}$. We choose one of them as the projection center $p$.  The resulting surface $S_p$ is the union of two surfaces $S_1$ and $S_2$. Both $S_1$ and $S_2$ are cones over over the same elliptic curve $C=S_1\cap S_2$.  By the Mayer-Vietoris sequence, we have 
$$\dots \to \sH^1(S_1)\oplus \sH^1(S_2)\to \sH^1(C)\to \sH^2(S_p)\to \sH^2(S_1)\oplus \sH^2(S_2)\to \dots$$
The resolution of $S_i$ (for $i=1,2$) is a ruled surface $\widetilde{S_i}$ over $C$. Moreover, the exceptional divisor of $\widetilde{S_i}\to S_i$ is isomorphic to $C$. It then follows that $\sH^1(S_i)=0$ and $\sH^2(S_i)$ is $1$-dimensional, carrying a pure Hodge structure of type $(1,1)$. In conclusion, we have $\Gr_1^W\sH^2(S_p)\cong \sH^1(C)$ and the claim follows.
\end{proof}

Similarly, for Type III  fourfolds we obtain:
\begin{lemma}\label{lemmacase3}
Let $X_0$ be a Type III fourfold with closed orbit. Then  $\Gr_2^W\sH^4(X_0)\neq 0$. 
\end{lemma}
\begin{proof}
By proposition \ref{propsingtype} (ii), except the case $\zeta$, we can find a singular point $p\in X_0$ not lying on a line of singularities. We can further assume that $p$ is of type $A_{\infty}$. As in  \ref{lemmacase2}, $S_p$ has insignificant limit singularities. Since some of the singularities of $S_p$ are degenerate cusps, $S_p$ is a Type III degeneration of $K3$ surface. It remains to consider the case $\zeta$. Assume $X_0$ is given by the equation (\ref{eqzeta}), and let $p=(1:0:\dots:0)$. The items (iii) and (iv) of proposition \ref{propblowup} do not hold for $S_p$ ($p$ lies on a line of singularities). In particular, since some of the singularities of $S_p$ are not hypersurface singularities,  we can not apply the results of Shah \cite{shahinsignificant} to get information about $\sH^2(S_p)$. On the other hand, since the singularities of $S_p$ are complete intersections, all the statements of  \ref{propblowup}, except (iii) and (iv), are still valid. It follows then that the relation  (\ref{eqreduction3}) between the Hodge structures of $X_0$ and $S_p$ also  holds. Finally, we find the Hodge structure on $\sH^2(S_p)$ by a direct computation as follows. The surface $S_p$ is the complete intersection $(x_4x_5=0,\ x_1x_2x_3=0)$ in $\bP^4$ (with homogeneous coordinates $(x_1:\dots:x_5)$). Let $S_1$ and $S_2$  be the  restrictions of $S_p$ to the hyperplanes ($x_4=0$) and ($x_5=0$) respectively. The surface $S_p$ is obtained by gluing $S_1$ and $S_2$ along the curve $C$ given by ($x_1x_2x_3=0$) in the plane ($x_4=x_5=0$). Thus, by Mayer--Vietoris, we have 
$$\dots \to \sH^1(S_1)\oplus \sH^1(S_2)\to \sH^1(C)\to \sH^2(S_p)\to \dots$$
Since $S_1$ and $S_2$ are normal crossing varieties, it follows easily that $\sH^1(S_i)=0$ for $i=1,2$ (e.g. \cite[pg. 111]{clemensschmid}). Also,  since $C$ is a cycle of rational curves, $\Gr_0^W\sH^1(C)$ is $1$-dimensional. Thus, $\Gr_0^W\sH^2(S_p)\neq 0$ and the lemma follows. 
\end{proof}

\begin{remark}\label{remresolution}
It is possible to compute the mixed Hodge structure of a Type II fourfold $X_0$ by means of a resolution of the singularities of $X_0$. This approach identifies the Hodge structure on $\Gr_3^W\sH^4(X_0)$ with the Hodge structure of an elliptic curve $C$ naturally associated to the singularities of $X_0$. For example, in the simplest case, the case $\phi$, $X_0$ is singular along an elliptic normal curve $C$. The singularities of $X_0$ are resolved by the blow-up of $C$. The resulting exceptional divisor $E$ is a quadric bundle over $C$. Since the intermediate Jacobian of $E$ is (up to isogeny) the Jacobian of $C$, we obtain an identification between $\Gr_3^W\sH^4(X_0)$ and $\sH^1(C)(-1)$. The situation is similar in the case when $X_0$ is singular along a rational normal curve with $4$ special points on it (see Prop. \ref{propsingtype} (i)). Finally, in the case of $\widetilde{E}_r$ singularities, a weighted blow-up of the $\widetilde{E}_r$ singularity produces a fourfold singular along an elliptic curve, bringing us to the situation discussed above.
 \end{remark}

%%%%%%%%%%%%%%%%%%%%%%%%%%%%%%%%%%%%%%%%%%%%%%%%%%%%%%%%%%%%%%
\section{Degenerations to Type IV fourfolds}\label{sectresolution}
Due to the presence of bad singularities for Type IV  fourfolds, the period map has singularities along the curve $\chi \subset \overline{\calM}$. In this section, we note that a certain blow-up $\widetilde{\calM}$  of $\overline{\calM}$ along $\chi$ essentially resolves these singularities. We say essentially because our result (see \ref{thmmonodromy2}) is weaker than an extension statement. However, a posteriori, one can indeed conclude that the blow-up $\widetilde{\calM}$ resolves the period map (see remark \ref{remresolution2}).

%More precisely, in \S\ref{sectkirwan}, using the fact that $\chi$ parametrizes fourfolds with large stabilizer we blow-up $\overline{\calM}$  along $\chi$ according to the canonical Kirwan desingularization procedure (see \cite{kirwan}). Next, in \S\ref{sectmodular}, we prove that the resulting space $\widetilde{\calM}$ has a modular interpretation. Essentially, the exceptional locus of $\widetilde{\calM}\to \overline{\calM}$ parametrizes quadric bundles canonically associated to degree two K3 surfaces. Additionally, the singularities allowed for these quadric bundles are of the same type as those of Type I--III fourfolds. Thus, via the blow-up procedure we replace the Type IV fourfolds by Type I--III fourfolds, for which we can control the monodromy as in section \ref{sectmonodromy}. 

\subsection{The Kirwan desingularization along $\chi$}\label{sectkirwan} We start by recalling the setup of Kirwan \cite{kirwan}. Let $G$ be a reductive group, acting on a smooth projective variety $P$. For $R$ a reductive subgroup of $G$ we denote by  $Z_R^{ss}$ the subset of semistable points stabilized by $R$. Then, to each conjugacy class of reductive subgroups with $Z_R^{ss}\neq \emptyset$ there is associated  a natural blow-up of the GIT quotient $P\gquot G$. The center of the blow-up is the closed subset $(G\cdot Z_R^{ss})\gquot G$.  The resulting space is again a GIT quotient, but with simpler singularities. By performing this procedure a finite number of times starting with the highest dimensional subgroup $R$, one obtains a space with only  quotient singularities.

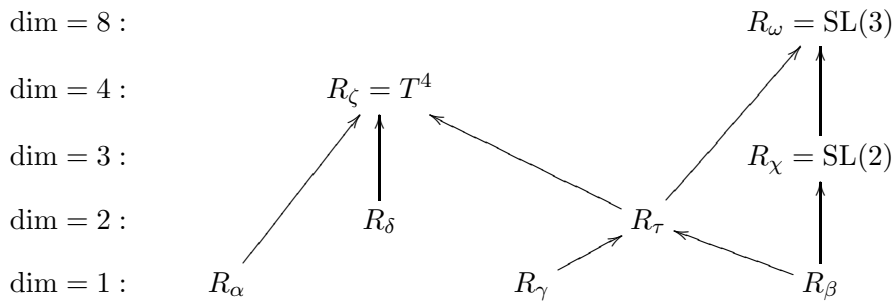
\begin{figure}[htb]
\begin{center}
$$
\xymatrix@R=.25cm{
{\dim=8:}&&&&&{R_{\omega}=\SL(3)}\ar@{<-}[dd] \ar@{<-}[dddl]\\
{\dim=4:}&&{R_{\zeta}=T^4}\ar@{<-}[dd] \ar@{<-}[dddl]\ar@{<-}[ddrr]\\
{\dim=3:}&&&&&{R_{\chi}}=\SL(2)\ar@{<-}[dd] \\
{\dim=2:}&&{R_{\delta}}&&{R_{\tau}}\ar@{<-}[dl]\ar@{<-}[dr]&\\
{\dim=1:}&{R_{\alpha}}&&{R_{\gamma}}&&{R_{\beta}}
}
$$
\end{center}
\caption{The adjacencies of connected reductive subgroups with $Z_R^{ss}\neq \emptyset$}\label{listr}
\end{figure}

\smallskip

In our situation $\overline{\calM}= \bP(\Sym^3 W)\gquot G$, where $G=\SL(6)$ and $W$ is the standard $G$-representation. The  conjugacy classes of subgroups $R$ are easily determined from the results of \cite{gitcubic}. Specifically, the minimal subgroups are $R_\alpha,\dots,R_\delta$  corresponding to the stabilizer of a generic point in the boundary strata $\alpha$--$\delta$ (see \cite[Sect. 4]{gitcubic}). These boundary strata further specialize to the strata $\tau$, $\chi$, $\zeta$, and $\omega$. The resulting subgroups (except when specified otherwise, these subgroups are tori of the indicated dimension) and  the corresponding inclusions are given in figure \ref{listr} (for a similar situation see \cite[pg. 66]{kirwanhyp}). 

\smallskip

We consider here only the first two steps of the Kirwan partial desingularization of $\overline{\calM}$. Namely, we blow-up $\overline{\calM}$ with respect  to $R_{\omega}$, followed by the blow-up with respect to $R_{\chi}$. We denote the resulting variety $\widetilde{\calM}$, and by $\widehat{\calM}$ the intermediary blow-up. Thus, $\widehat{\calM}$ is the blow-up of $\overline{\calM}$  in the point $\omega$, and $\widetilde{\calM}$ is the blow-up of  $\widehat{\calM}$  along the strict transform of $\chi$.

\subsubsection{The blow-up of the point $\omega$}\label{sectblowup1}
The moduli of cubic fourfolds is obtained as the quotient of  the action of $G=\SL(6)$ on $P:=\bP(H^0(\bP^5,\calO_{\bP^5}(3)))=\bP(\Sym^3 W)$. Let $\omega\in \calM$ be the special boundary point. The preimage of $\omega$ in $P$ contains a unique closed orbit, say the orbit of $x\in P^{ss}$. 

\smallskip

The local structure (in the \'etale topology) of the  quotient $\overline{\calM}$ is described by Luna's slice theorem \cite[App. D]{GIT}. Specifically, there exists a $G_x$-invariant slice $W$ to the orbit $G\cdot x$. The slice $W$ can be taken to be a smooth, affine, locally closed subvariety of $P$ such that $U=G\cdot W$ is open in $P$. We then have the following commutative diagram with Cartesian squares: 
$$
\begin{CD}
G\times_{G_x}\calN_x@<\text{\'etale}<<G\times_{G_x}W@>\text{\'etale}>>U@.\subset@.P\\
@VVV                               @VVV@VVV@.@VVV\\
\calN_x\gquot G_x@<\text{\'etale}<<W\gquot G_x@>\text{\'etale}>>\ U\gquot G\ @.\subset@.\ P\gquot G\ @.\cong@. \overline{\calM}\end{CD}
$$
where $G_x$ denotes the stabilizer of  $x$, and $\calN_x$ the fiber of the  normal bundle of the orbit $G\cdot x$. The Kirwan desingularization is compatible with Luna's slice. In particular, the exceptional divisor of the blow-up of $\overline{\calM}$ at the point corresponding to  $G\cdot x$ is the quotient $\bP(\calN_x)\gquot G_x$ (\cite[Rem. 6.4]{kirwan}). 

\smallskip

We recall that the point $\omega$ corresponds to  the secant variety $X$ of the Veronese surface $S$ in $\bP^5$. It follows that $G_x^0\cong \SL(3)$,  acting in the standard way on the singular locus $S$ of $X$. The Veronese surface $S\subset \bP^5$  is the image of $\bP^2$ embedded by $\calO_{\bP^2}(2)$ in $\bP^5$. Thus, there is a natural identification  $W\cong \Sym^2 V$ of $\SL(3)$-representations such that $S\subset \bP(W)$ is left invariant (where $V$ is the standard $\SL(3)$ representation). We obtain the decomposition in irreducible summands 
\begin{equation}\label{decomposition}
\Sym^3 W\cong \Sym^3 \Sym^2 V\cong  \Sym^6 V\oplus \Gamma_{2,2}\oplus \bC
\end{equation}
(cf. \cite[(13.15)]{fultonharris}). It then follows easily that the natural representation of $\SL(3)$ on the normal slice $\calN_x$ is isomorphic to $\Sym^6 V$.   In conclusion, from Luna's slice theorem and this computation we have obtained the following precise description of the blow-up of the point $\omega$:

\begin{corollary}\label{cormodular}
The point $\omega\in \overline{\calM}$ has an  \'etale (or analytic) neighborhood  isomorphic to the affine cone over the moduli space of plane sextic curves. The Kirwan blow-up of $\overline{\calM}$ at the point $\omega$ is isomorphic when restricted to this neighborhood to the natural blow-up of the vertex of the cone.  In particular the exceptional divisor is isomorphic to the GIT quotient for plane sextics. 
\end{corollary}

\subsubsection{The blow-up of the curve $\chi$}\label{sectblowup2} We repeat the computation from the previous section for the strict transform of the curve $\chi$. We recall that the singular locus of a fourfold giving a minimal orbit of type $\chi$ is a rational normal curve of degree $4$. It follows that (at least generically) $G_x\cong \SL(2)$, for $x\in P$ with closed orbit and mapping to $\chi$. Similarly to the case $\omega$ we obtain:

\begin{lemma}
Let $x\in P$ be a semistable point with closed orbit, mapping to $\chi$.  Then the natural representation of $G_x^0\cong \SL(2)$ on the normal  $\calN_x$ is isomorphic to 
\begin{equation}\label{repnormal}
\Sym^{12} V\oplus \Sym^8 \oplus \bC,
\end{equation}
where $V$ denotes the standard $\SL(2)$ representation.
\end{lemma}
\begin{proof}
A fourfold of type $\chi$ is singular along a rational normal curve of degree $4$ contained in a hyperplane of $\bP^5$. It follows (as in \S\ref{sectblowup1}) that we can naturally identify $W$ to the $\SL(2)$ representation $\Sym^4 V\oplus \bC$. We then have 
\begin{eqnarray}
\Sym^3 W&\cong& \Sym^3 (\Sym^4 V\oplus \bC)\cong \Sym^3 (\Sym^4 V)\oplus \Sym^2 (\Sym^4 V) \oplus \Sym^4 V \oplus \bC\\
\nonumber &\cong&\Sym^{12} V\oplus (\Sym^8 V)^{\oplus 2}\oplus \Sym^6 V\oplus (\Sym^4 V)^{\oplus 3} \oplus \bC^{\oplus 3}
\end{eqnarray}
The lemma now follows from the identification of the summands in the normal  sequence:
$$0\to \calT_{G\cdot x, x}\to \calT_{P,x}\to \calN_x\to 0$$
(e.g. the representation on $\calT_{G\cdot x, x}$ is computed by using $G\cdot x\cong G/G_x$). 
\end{proof}

To interpret geometrically the result of the previous lemma, we recall the situation of plane sextics studied by Shah \cite{shah}. Namely, to resolve the rational map from the GIT quotient of plane sextics to the period domain of degree two K3 surfaces, one needs to blow-up the point corresponding to the triple conic. The stabilizer is again $\SL(2)$, and the representation on the normal slice is 
\begin{equation}\label{repnormal2}
\Sym^{12} V\oplus \Sym^8 V
\end{equation}
The two summands of (\ref{repnormal2}) correspond in Shah's notation (\cite[pg. 498-500]{shah})  to  $\Phi$ and $\Xi$ respectively. In conclusion, by comparing (\ref{repnormal}) and  (\ref{repnormal2}), we can say that the blow-up of $\chi$  is locally the product of the affine line with the exceptional divisor $\fN$ obtained by Shah \cite[\S5]{shah}. The divisor $\fN$ corresponds to elliptic degree two K3 surfaces. These special degree two K3 surfaces are not double covers of $\bP^2$, but instead are double covers of the Hirzebruch surface $\bF_4$ embedded in $\bP^5$ as the cone $\Sigma_4^0$ over a rational normal curve of degree $4$. 

\begin{remark}\label{remunigonal}
It is well known that for a degree two K3 surface $(S,H)$ the polarization $2H$ is base point free and defines a morphism to $\bP^5$. In the generic case the image of $S$ is the Veronese surface, and in the special (unigonal) case the image is $\Sigma_4^0$.  In both cases the branch curve is cut by a cubic fourfold. Thus, from our point of view the two cases behave similarly.
\end{remark}

%%%%%%%%%%%%%%%%%%%%%%%%%%%%%%%%%%%%%%%%%%%%%%%%%%%%%%%%%%%%%%
\subsection{Modular interpretation of the blow-up $\widetilde{\calM}\to \overline{\calM}$}\label{sectmodular}
In \S\ref{sectblowup1} and \S\ref{sectblowup2} we have noted the natural occurrence of the moduli space of degree two K3 surfaces in the blow-up  of $\overline{\calM}$ along $\chi$. The statement there is purely representation theoretic. We now interpret geometrically this observation and relate it to the Hodge theory of the degenerations to Type IV fourfolds. Most of the discussion here is about the case $\omega$.  Since the case $\chi$ is quite similar (see \S\ref{sectblowup2}, especially remark \ref{remunigonal}),  there is no difficulty in adapting the arguments to the case $\chi$. In particular we note that the semi-stable reduction in this case is essentially contained in Allcock--Carlson--Toledo \cite[Sect. 5]{allcock3fold}.

\smallskip

 To start, we recall the computation of the semi-stable reduction for a general pencil $\calX\to \Delta$  of cubic fourfolds degenerating to the secant variety $X_0$ of  Veronese surface $V$ (cf. Hassett  \cite[\S4.4]{hassett}). The total space $\calX$ is singular along the degree $12$ curve $C=X_t\cap V$ (a plane sextic). After a base change of order $2$, and the blow-up of the family $\calX$ along the Veronese surface $V$, we obtain a semi-stable model $\calX'\to \Delta$. The new central fiber $X_0'$ is the union of two smooth components $\overline{X}_0$ and $E$, intersecting transversely in  a $\bP^1$-bundle $E_0$ over $\bP^2$. The fourfold $\overline{X}_0\cong \textrm{Hilb}^2(\bP^2)$ is the natural resolution of $X_0\cong \Sym^2(\bP^2)$. The other component $E$ of the central fiber is  a quadric bundle over $V\cong \bP^2$ with discriminant the curve $C$. As it is standard, to $E$ one associates the degree two $K3$ surface $S$ obtained as the double cover of $\bP^2$ branched along $C$ (intrinsically $S\xrightarrow{2:1}\bP^2$ parameterizes the rulings of the fibers of $E\to \bP^2$). The Hodge structure of $E$  (and that of $X_0'$) is essentially the Hodge structure of $S$. Specifically, in the generic case, the Hodge structure on the transcendental sublattice of $\sH^4(E,\bZ)$ is a Tate twist of the Hodge structure on the primitive cohomology of $S$. It follows then that $\sH^2_0(S)(-1)$ embeds in the limit Hodge structure $\sH^4_{\lim}$. This implies that $\sH^4_{\lim}$ is pure and the corresponding period point belongs to the divisor $\calH_{\infty}/\Gamma$. By interpreting 
 the pencil $\calX\to \Delta$ as a generic arc meeting the exceptional divisor (corresponding to $\omega$) of the blow-up $\widehat{\calM}\to \overline{\calM}$, we conclude that the monodromy around this exceptional divisor is of order $2$ and corresponds to a reflection in a long root (see remark \ref{remloop}).

 \smallskip
 
 The above computation still remains valid for any arc intersecting transversally the exceptional divisor of the blow-up of $\omega$ (N.B. the computations should be done in \'etale slice modeling the local structure of $\overline{\calM}$ near $\omega$). Thus, we can interpret this exceptional divisor as parametrizing fourfolds of type $X_0'=\overline{X}_0\cup_{E_0} E$, where $E$ is a (possibly singular) quadric bundle over $\bP^2$ (i.e. $E\to \bP^2$ is an algebraic  fiber space with generic fiber a smooth quadric in $\bP^3$). This identification is rather tautological. Namely, by \S\ref{sectblowup1}, the points of the exceptional divisor of $\widetilde{\calM}\to \overline{\calM}$ naturally correspond to isomorphism classes of plane sextics. The following lemma shows that to a plane sextic we can associate a quadric bundle $E$ in a compatible way with our identification.  
 
\begin{lemma}\label{lemmaquadric} 
 The quadric bundles $E$ occurring in the semistable reduction of degenerations to the secant to Veronese are given by sections $q\in H$, where 
 $$H:=\sH^0\left(\bP^2,\Sym^2(\calN^*_{V/\bP^5}\oplus \calO_{\bP^2})\otimes \calO_{\bP^2}(6)\right).$$ 
 Furthermore, with respect to the natural $G=\SL(3)$ action on $H$, we have the decomposition:
\begin{equation}
H\cong \sH^0\left(\bP^2,\Sym^2(\calN^*_{V/\bP^5})\otimes \calO_{\bP^2}(6)\right) \oplus \sH^0\left(\bP^2,\calN_{V/\bP^5}\right)\oplus \sH^0\left(\bP^2,\calO_{\bP^2}(6)\right). 
\end{equation}
In particular, $q$ can be taken of type $q=q_0+f$ with $q_0\in \sH^0\left(\bP^2,\Sym^2(\calN^*_{V/\bP^5})\otimes \calO_{\bP^2}(6)\right)$ a tautological $G$-invariant section and $f\in \sH^0\left(\bP^2,\calO_{\bP^2}(6)\right)$ the equation of the  discriminant curve $C$. 
\end{lemma}
\begin{proof}
 As the semi-stable reduction is obtained by blowing-up the Veronese surface,  we have  $E\subset\bP(\calN^*_{V/\bP^5}\oplus \calO_{V})$ (see also \cite[Lemma 4.4.3]{hassett}). The quadric bundle $E$ corresponds to a  section $q\in \sH^0\left(\bP^2,\Sym^2(\calN^*_{V/\bP^5}\oplus \calO_{V})\otimes \calO_{\bP^2}(\alpha)\right)$ for some $\alpha\in \bZ$ (see \cite[Prop. 1.2]{beauvilleprym}). Then, the discriminant curve is  given by the induced section:  
$$\Delta\in \sH^0\left(\bP^2,\Sym^2\left(\det\left(\calN^*_{V/\bP^5}\oplus \calO_{V}\right)\right)\otimes \calO_{\bP^2}(4\alpha)\right).$$

The following basic facts about the normal (and conormal) bundle of the Veronese surface $V$ in $\bP^5$  are well known (e.g. \cite[Lemma 2.7]{einsb}):
\begin{itemize}
\item[i)] $\calN_{V/\bP^5}^*\cong \Sym^2(\Omega^1_{\bP^2})$; % (cf. \cite[Remark 2]{egps});
\item[ii)] $\calN_{V/\bP^5}^*(2)\cong \calN_{V/\bP^5}(-1)$ (the twist is taken with respect to $\calO_{\bP^5}(1)$);
\item[iii)] The isomorphism of ii) induces a $G$-invariant section $q_0$ via the natural inclusion:
$$q_0\in \sH^0\left(\Sym^2(\calN_{V/\bP^5}(-1))\otimes \calO_V(-1)\right)\hookrightarrow \Hom\left(\calN_{V/\bP^5}^*(1), \calN_{V/\bP^5}(-1)\otimes \calO_V(-1)\right)$$ 
\end{itemize}
Using ii), we rewrite iii) as $q_0\in  \sH^0(\Sym^2(\calN^*_{V/\bP^5})\otimes \calO_{\bP^2}(6))$. The lemma now follows from \cite[Lemma 4.4.3]{hassett}  by noting the surjection $\sH^0(\bP^5,\calO_{\bP^5}(3))\twoheadrightarrow \sH^0(V,\calO_{V}(3))\cong \sH^0(\bP^2,\calO_{\bP^2}(6))$.
\end{proof}

The previous lemma associates in a  canonical way to a plane sextic, a quadric bundle over $\bP^2$. Since we are considering also singular sextics, we need to understand the possible singularities of the resulting quadric bundle. 

\begin{lemma}\label{lemmasing}
Let $E$ be a quadric bundle as in lemma \ref{lemmaquadric}. Then the singularities of $E$ are in one-to-one correspondence, including the analytic type, with the singularities of the discriminant curve $C$.
\end{lemma}
\begin{proof}
We choose $U\subset \bP^2$ an open affine over which  $\calN_{V/\bP^5}$ is trivialized. Locally over $U$, $E_0\subset  \bP^2\times U\to U$ is given by $x_1^2+x_0x_2$ and $E\subset \bP^3\times U\to U$  by $f(x,y)t^2+x_1^2+x_0x_2$, where $(x,y)$ are affine coordinates on $U$, $f$ is the equation of the discriminant curve, and $(x_0:x_1:x_2:t)$ are appropriate homogeneous coordinates on $\bP^3$. Since $E_0$ is the smooth hyperplane section ($t=0$) of $E$, we can assume $t=1$. Clearly, the singularities of $E$ correspond to the singularities of $C$.
\end{proof}

In the case of arcs $\tau:\Delta\to \widehat{\calM}$ meeting non-transversally the exceptional divisor of the blow-up of the point $\omega\in\overline{\calM}$, one might need several blow-ups to obtain the quadric bundle $E$. More precisely, starting with a $1$-parameter degeneration $\calX\to \Delta$ to $\omega$, repeating several times the computation from the transversal case, we arrive to a partial semi-stable reduction $\calX'\to \Delta$ such that one of the components  of the central fiber is a quadric bundle $E$  attached to the other components along a smooth divisor $E_0$,  and such that the condition of semi-stable family for $\calX'$ holds away from $E$. As before, the quadric bundle $E$ is determined by the point where the arc meets the exceptional divisor. By standard arguments (e.g. \cite[Prop. 2.1]{shah}), we can further assume that discriminant curve $C$ associated to $E$ is  a semi-stable plane sextic with closed orbit (see Shah \cite[Thm. 2.4]{shah}). Then, depending on the singularities of $E$ (or equivalently $C$) we call the central fiber $X_0'$ of $\calX'$ a fourfold of Type I'--IV', where the  types are defined by:
\begin{itemize}
 \item[(I')]  the cases when $E$ has at worst simple singularities; 
 \item[(II')] the cases when $E$ has at worst $\widetilde{E_r}$, $A_{\infty}$ or $D_{\infty}$ singularities (corresponding to the Type II in Shah \cite[Thm. 2.4]{shah});
 \item[(III')] the cases when $E$ has  at worst degenerate cusp singularities (Type III in \cite{shah}); 
 \item[(IV')] the cases corresponding to the intersection of the exceptional divisor with the strict transform of the curve $\chi$ (i.e. the discriminant sextic is the triple conic).
\end{itemize}

We conclude that any degenerations of cubic fourfolds to the secant variety of the Veronese surface can be filled in with a fourfold of type I'--IV' (as defined above). Since the interesting cohomology of $X_0'$ comes from $E$, which in  turn comes  from the cohomology of a degeneration of $K3$ surfaces with insignificant singularities, we can control the monodromy for degenerations with central fiber of Type I'--III' as in section \ref{sectmonodromy} (i.e. the index of nilpotency is the numerical label of the type). By considering the blow-up of $\widetilde{\calM}\to \widehat{\calM}$ of  the strict transform of $\chi$ in $\overline{\calM}$, we can eliminate the remaining cases of Type IV (corresponding to $\chi\setminus \omega$) and  the Type IV', introducing instead fourfolds analogous to the Types I'--III' above  (this is essentially Shah \cite[Thm. 4.3]{shah}). In conclusion, using the blow-up $\widetilde{\calM}\to \overline{\calM}$, any $1$-parameter family of smooth cubic fourfolds degenerating to a Type IV fourfold can be filled in with a fourfold of Type I'--III', for which we can control the monodromy. In particular, for the Type II' and III' cases the limit period belongs to the boundary of $(\calD/\Gamma)^*$. On the other hand, as already noted above in the generic case, the limit period for a Type I' fourfold belongs to $\calH_{\infty}/\Gamma$.

%%%%%%%%%%%%%%%%%%%%%%%%%%%%%%%%%%%%%%%%%%%%%%%%%%%%%%%%%%%%%%
\begin{theorem}\label{thmmonodromy2}
Let $f:\calX\to \Delta$ be a $1$-parameter degeneration with finite monodromy of cubic fourfolds such that the central fiber is of Type IV. Then, the period point corresponding to the limit  Hodge structure $\sH^4_{\lim}$ (which is automatically pure) belongs to $\calH_{\infty}/\Gamma$. 
\end{theorem}
\begin{proof}
 For simplicity, we restrict the discussion to the case $\omega$; the computation in the case $\chi$ is similar (in fact, $\omega$ is the degenerate case of $\chi$; one obtains an additional class in the case $\chi$). Assume that the central fiber  $X_0$ of the degeneration is the secant variety of  the Veronese surface $S\hookrightarrow \bP^5$. We have $X_0\cong \Sym^2(\bP^2)$. Since  $X_0$ is the quotient of $\bP^2\times \bP^2$ by the finite group $\mu_2$, the cohomology (over $\bQ$) of $X_0$  is the $\mu_2$-invariant cohomology of $\bP^2\times \bP^2$. It follows that $\sH^4(X_0)$ is $2$-dimensional carrying a pure Hodge structure of type $(2,2)$. There are two natural generators $h$ and $e$ for $\sH^4(X_0)$: the square of the hyperplane section and the class corresponding to the  diagonal in $\bP^2\times \bP^2$. Via the natural double cover $\bP^2\times \bP^2\to X_0$, one computes the intersection matrix of $h$ and $e$ to be $\left(\begin{array}{cc} 3&2\\
2&2\end{array}\right)$ (since $X_0$ is the quotient of a smooth variety by a finite group, the intersection makes sense). 
By the discussion of section \ref{sectperiod} (esp. lemma \ref{lemmaloop}), the limit period point belongs to $\calH_{\infty}/\Gamma$ provided that the algebraic classes $h$ and $e$ are preserved in the limit $\sH^4_{\lim}$. 

\smallskip

Let $\calX\to \Delta$ be a  degeneration of cubic fourfolds to $X_0$, the secant of the Veronese surface $S$. Assume that $\calX'\to \Delta$ is a semi-stable reduction of the degeneration. We can write the new central fiber as $X_0'=\overline{X}_0\cup E$, where $\overline{X}_0$ is the component dominating $X_0$ (thus a resolution of $X_0$) and $E$ the union of the other components in the central fiber $X_0'$. Let $E_0=\overline{X}_0\cap E$.  We get a  diagram:
\begin{equation}\label{eqres2}
\xymatrix{
E\ar@{->}[dr]&E_0\ar@{->}[d]\ar@{_{(}->}[l]\ar@{^{(}->}[r]&\overline{X}_0\ar@{->}[d]\\
&S\ar@{->}[r]&X_0
}
\end{equation} 
%Since we are only interested in the contribution of  $\sH^n(X_0)$ to the top weight part of the mixed Hodge structure on $\sH^n(X_0')$, there is no loss of generality in working as if there are only two components $\overline{X}_0$ and $E$. In general, $E$  is a union of normal crossing components $X_1\cup\dots \cup X_k$ (and similarly for $E_0$),.
The limit Hodge structure $\sH^4_{\lim}$ is computed by the Clemens--Schmid exact sequence  for the semi-stable family $\calX'\to \Delta$. Specifically, under the assumption of finite monodromy, we have: 
\begin{equation}\label{eqcs2}
0\to \sH^2_{\lim}\to \Gr^W_{-6}\ \sH_6(X_0')\xrightarrow{\alpha}\Gr^W_4  \sH^4(X_0')\to \sH^4_{\lim}\to 0
\end{equation}
(\cite[pg. 109]{clemensschmid}). The relevant part of the cohomology of $X_0'$ is computed by the exact sequence: 
\begin{equation}\label{mv5}
0\to \Gr^W_4 \sH^4(X_0')\to   \sH^4(\overline{X}_0)\oplus \Gr^W_4 \sH^4(E)\to\Gr^W_4  \sH^4(E_0)
\end{equation}
 Combining (\ref{mv5}) with the exact sequence (\ref{mayervietoris}) for the resolution $\overline{X}_0\to X_0$ (see diagram (\ref{eqres2})), we get a natural injective map $\sH^4(X_0)\hookrightarrow  \Gr^W_4 \sH^4(X_0')$. Since the morphism $\alpha: \sH_6(X_0')\to   \sH^4(X_0')$ is obtained as the composition of the Poincare duality map $ \sH_6(X_0')\cong \sH_6(\calX')\to \sH^4(\calX',\partial \calX')$ with the restriction map  $\sH^4(\calX',\partial \calX')\to \sH^4(\calX')\cong \sH^4(X_0')$ (\cite[pg. 108]{clemensschmid}), no class $s_0\in \sH^4(X_0')$ that can be extended to a non-zero class $s_t\in \sH^4(X_t')$ belongs to the image of $\alpha$. The classes coming from the embedding $\sH^4(X_0)\hookrightarrow \sH^4(X_0')$ extend to the nearby $X_t'$ fibers for the simple reason that $\sH^4(X_0)$ is generated by the polarization class and a class supported on the smooth locus of $X_0$. Thus, the image of $\sH^4(X_0)$ in $\sH^4_{\lim}$ is two dimensional and the claim follows.
 \end{proof}

\section{Proof of the main results}\label{sectproofs}
\subsection{Proof of Theorem \ref{mainthm1}} 
The period map $\calP_0:\calM_0\to \calD/\Gamma$ is a birational morphism of quasi-projective varieties with the image contained in $(\calD\setminus (\calH_{\Delta}\cup \calH_{\infty}))/\Gamma$. By proposition \ref{extension1},  $\calP_0$ extends to  a regular morphism $\calP:\calM\to \calD/\Gamma$ over the simple singularities locus $\calM\subset\widetilde{\calM}$. Additionally, the image of the discriminant is included in $\calH_{\Delta}$. Now, theorem \ref{mainthm1} follows provided the properness of the extended period map $\calP:\calM\to (\calD\setminus \calH_\infty)/\Gamma$. 

By the valuative criterion, it suffices to consider $1$-parameter families $\calX^*\to \Delta^*$  of smooth cubic fourfolds over the punctured disk $\Delta^*$ with the limit period point in $(\calD\setminus \calH_\infty)/\Gamma$. The properness means that $\calX^*\to \Delta^*$ can be filled in (after a possible finite base change) to a family $\calX\to \Delta$ of cubic fourfolds with at worst simple singularities. By using the GIT quotient $\overline{\calM}$, we can extend $\calX^*\to \Delta^*$ to a family $\calX\to \Delta$ of cubic fourfolds such that the central fiber $X_0$  is a semi-stable cubic fourfold with closed orbit (see \cite[Prop. 2.1]{shah}). If $X_0$ is of Type I (see \S\ref{gitresults}), we are done. Thus, assume $X_0$ is of Type II--IV. Since the limit period point for  $\calX^*\to \Delta^*$ belongs to the period domain, it is known that the monodromy of the family $\calX^*$ is finite (see \cite[Thm. 13.4.5]{carlsonperiod}). By theorem \ref{thmmonodromy}, we conclude that $X_0$ can not be of Type II or III.  Similarly, the Type IV case is excluded by theorem \ref{thmmonodromy2} (i.e. the limit period belongs to $\calH_{\infty}/\Gamma$ in this case). \qed

\subsection{Proof of the Theorem \ref{mainthm2}} The statement follows directly from theorem \ref{mainthm1} and the general results of Looijenga \cite[Thm. 7.6]{looijengacompact}. In fact, the situation is formally the same as that of low degree K3 surfaces, case that was analyzed in detail by Looijenga \cite[\S8]{looijengacompact} (esp. \cite[Thm. 8.6]{looijengacompact}). A sketch of the argument is given below.

\smallskip

Let $G=\SL(6)$, and $P=\bP \Sym^3 W$ be the linear system of cubics in $\bP^5$. On $P$ we consider the natural polarization $\eta=\calO_P(n)$ (for some $n$ to be chosen later), and define $U\subset P$ to be the $G$-invariant open subset parametrizing cubic fourfolds with at worst simple isolated singularities. The GIT results of \cite{gitcubic} imply in particular that:
\begin{itemize}
\item[i)]  $P\setminus U$ has high codimension (at least $2$) in $P$ ; 
\item[ii)] the points of $U$ are stable with respect to the action of $G$ on $P$ (\cite[Thm. 1.1]{gitcubic}).
\end{itemize}
From the general theory, we note that the GIT quotient $\overline{\calM}$ is isomorphic to the $\Proj$ of the ring of $G$-invariant sections (i.e. $\Proj(\oplus_{k\ge 0} H^0(P,\eta^{\otimes k})^G)$) and that some power of $\eta$ descends to an ample line bundle on $\overline{\calM}$.

\smallskip

On the $\calD/\Gamma$ side, there exists a natural automorphic line bundle $\bL$ on the domain $\calD$, which descends to an ample line bundle $\calL$ on $\calD/\Gamma$ (by the Baily--Borel theory). Let $\calD^0=\calD\setminus \calH_{\infty}$ and $X^0=\calD^0/\Gamma$. The content of theorem \ref{mainthm1} is that the period map induces an isomorphism between the quasi-projective varieties $\calM=U/G$ and $X^0=\calD^0/\Gamma$. A basic observation of Looijenga is that, for hypersurfaces, one has  a stronger statement: the isomorphism given by the period map is an isomorphism of polarized varieties, i.e. the period map gives an  identification of $(U,\eta_{\mid U})/G$ with $(X^0,\calL_{\mid X^0})$.  Namely, given the identification between $\calM$ and $X^0$, the identification of line bundles follows from Griffiths residue theory: the period point corresponding to a cubic fourfold is determined by the $1$-dimensional subspace $\sH^{3,1}$ of the middle cohomology, which in turn is generated by the residue of a rational form on $\bP^5$, form that depends only on the equation of the cubic (in particular, since this rational form depends quadratically on the equation of the cubic, the right choice of $n$ in the definition of $\eta$ is $2$).

\smallskip

We conclude that the assumptions (in particular (i) and (ii) from above) of  \cite[Thm. 7.6]{looijengacompact} are satisfied. Thus, we  obtain an isomorphism of graded algebras:
$$\oplus_{k\ge 0} H^0(P,\eta^{\otimes k})^G\cong \oplus_{k\in \bZ} H^0(\calD^0,\calO(\bL)^{\otimes k})^\Gamma.$$
Taking $\Proj$ of both sides, one obtains an isomorphism between the GIT compactification $\overline{\calM}$ and the Looijenga compactification associated to $\calH_\infty$ (for a geometric description of this compactification see the next section, esp. diagram \ref{birationaldiag}). Thus, theorem \ref{mainthm2} holds.
 \qed

%%%%%%%%%%%%%%%%%%%%%%%%%%%%%%%%%%%%%%%%%%%%%%%%%%%%%%%%%%%%%%
\section{Arithmetic properties of the arrangement of hyperplanes $\calH_\infty$}\label{sectapplication}
We close by noting several arithmetic properties of the hyperplane arrangement $\calH_{\infty}$, which give (a posteriori) some information about  the GIT compactification $\overline{\calM}$. Namely, we recall that as a consequence of Theorem \ref{mainthm2} and of the general results of Looijenga \cite{looijengacompact} the birational map between the GIT compactification $\overline{\calM}$ and the Baily--Borel compactification $(\calD/\Gamma)^*$ fits into a diagram:
\begin{equation}\label{birationaldiag}
\xymatrix@R=.25cm{
{\widetilde{\calM}}\ar@{->}[dd]\ar@{->}[r] &\ar@{->}[dd]  \widehat{\calD/\Gamma}\\
\\
\overline{\calM}\ar@{-->}[r] &(\calD/\Gamma)^*
}
\end{equation}
The space $\widehat{\calD/\Gamma}$ is a small modification of  $(\calD/\Gamma)^*$ such that the Weil divisor defined by $\calH_{\infty}$ becomes $\bQ$-Cartier. Then $\widetilde{\calM}$  is obtained from $\widehat{\calD/\Gamma}$  by blowing-up the strict transform of the arrangement $\calH_{\infty}$ in the usual way (starting with the highest codimension self-intersection stratum). Finally, the morphism $\widetilde{\calM}\to\overline{\calM}$ contracts the hyperplane arrangement in the opposite direction. The net effect of this birational modification is to flip a codimension $k$ intersection of hyperplanes in $(\calD/\Gamma)^*$ to a ($k-1$)-dimensional stratum in $\overline{\calM}$. Since all the steps involved in passing from $(\calD/\Gamma)^*$ to $\overline{\calM}$  are explicit and of relatively simple arithmetic nature, we can recover some of the GIT results obtained in \cite{gitcubic} (see the discussion from \cite[Remark 3.2]{looijengaswierstra} for the similar situation of cubic threefolds). For instance,  the dimensions of the boundary strata from table \ref{tableboundary} can be obtained by understanding how the arrangement of hyperplanes $\calH_\infty$ meets the boundary of $(\calD/\Gamma)^*$.

\begin{remark}\label{remresolution2}
An unpublished argument of Looijenga (based on the universal property of the blow-up of the arrangement $\calH_{\infty}$, see \cite[Prop. 7.2]{looijengacompact}) shows that an immediate consequence of theorem \ref{mainthm2} is that the space $\widetilde{\calM}$ of diagram (\ref{birationaldiag}) is isomorphic to the Kirwan blow-up considered in section \ref{sectresolution}. This explains our choice of notation. However, in what follows, we do not need to assume this.  
\end{remark}

\subsection{The Baily--Borel Compactification of $\calD/\Gamma$}\label{sectbb}
The Baily--Borel theory compactifies the locally symmetric space $\calD/\Gamma$ to the normal projective variety $(\calD/\Gamma)^*$. The compactification is done by adding points, {\it Type III boundary components}, and curves, {\it Type II boundary components}. The Type II and III boundary components correspond $1$-to-$1$ to  the classes (modulo $\Gamma$) of the isotropic sublattices of $\Lambda_0$ of rank $2$ and $1$ respectively. Furthermore, the incidences between the boundary components correspond to the inclusions of the isotropic lattices. In our situation, we obtain:

\begin{proposition}\label{thmbbcompact}
The period domain $\calD/\Gamma$ for cubic fourfolds is compactified by adding six type II boundary components meeting in a single point, the unique type III component. The type II boundary components are distinguished by the root sublattice $R$ contained in  $E^\perp/E$, where $E$ is the rank $2$ isotropic sublattices associated to the component. The six possibilities for $R$ are $E_8^{\oplus 2}\oplus A_2$, $D_{16}\oplus A_2$, $E_7\oplus D_{10}$, $A_{17}$, $E_6^{\oplus 3}$ and $A_{11}\oplus D_7$. 
\end{proposition}
\begin{proof}
We first note that up to the action of $\Gamma$ there exists only one class of isotropic rank $1$ sublattices in $\Lambda_0$. 
Since  $\Lambda_0$ contains two hyperbolic planes, the claim follows from general lattice theory (e.g. \cite[Lemma 4.1.2]{scattone}). Thus, there exists a unique type III component. 

For the classification of rank $2$ isotropic sublattices in $\Lambda_0$ we  apply the approach of Scattone \cite{scattone}. Namely, let $E$ be an isotropic rank $2$ lattice. A standard invariant of $E$ is the  isometry class of the lattice $N:=E^\perp/E$. The lattice $N$ is positive definite of rank $18$ and in the same genus as $E_8^{\oplus 2}\oplus A_2$. By gluing $N$ and $E_6$ along the discriminant groups, we obtain an even unimodular lattice $L$ of rank $24$. Conversely, any lattice $N$ in the genus of $E_8^{\oplus 2}\oplus A_2$ can be realized  as the orthogonal complement of an embedding of $E_6$ into an even  unimodular lattice $L$ of rank $24$ (see \cite[Prop. 6.1.1]{scattone}). The possibilities for $L$ were classified by Niemeier:  there are $24$ isometry classes of even unimodular rank $24$ lattices,  which are distinguished by the sublattice 
spanned by the roots. It is clear that  $E_6$ can be embedded only into those unimodular lattices $L$ for  which the sublattice spanned by the roots  contains an $E_r$ summand (for $r=6,7,8$).  We get six possibilities for $L$ corresponding to the root lattices: $E_6^{\oplus 4}$, $E_6\oplus A_{11}\oplus D_7$, $E_7^{\oplus 2}\oplus D_{10}$, $E_8^{\oplus 3}$ and $E_8\oplus D_{16}$ (see \cite[\S18.4]{conway}).  Since $E_6$ embeds into an $E_r$ summand uniquely up to isometries, we obtain only six possibilities for $N$. We label 
(unambiguously)  the six cases by the root sublattice contained in $N$, i.e. $E_8^{\oplus 2}\oplus A_2$, $D_{16}\oplus A_2$, $E_7\oplus D_{10}$, $A_{17}$, $E_6^{\oplus 3}$ and $A_{11}\oplus D_7$ respectively.

To conclude the proof, we have to show that the classification of the possibilities for $N$ is equivalent to the classification of $E$. One direction is easy. Given $N$  in the genus of $E_8^2\oplus A_2$, there exists an isotropic lattice $E$ such that $N\cong E^\perp/E$. Namely, since $\Lambda_0$ is unique in its genus (see \cite{nikulin}), we have  $\Lambda_0\cong N\oplus U^2$ and the claim follows.  Conversely, $E$ is uniquely determined by $N$ as in Scattone \cite[Thm. 5.0.2 (1)]{scattone}. 
Finally, the classification of the isotropic lattices $E$ given above is (a priori) only up to the orthogonal group $\mathrm{O}(\Lambda_0)$. Since $\mathrm{O}^+(\Lambda_0)/\mathrm{O}^*(\Lambda_0)\cong \{\pm \mathrm{id}\}$, we obtain in fact the same classification modulo  $\Gamma\cong \mathrm{O}^*(\Lambda_0)$.
\end{proof}

  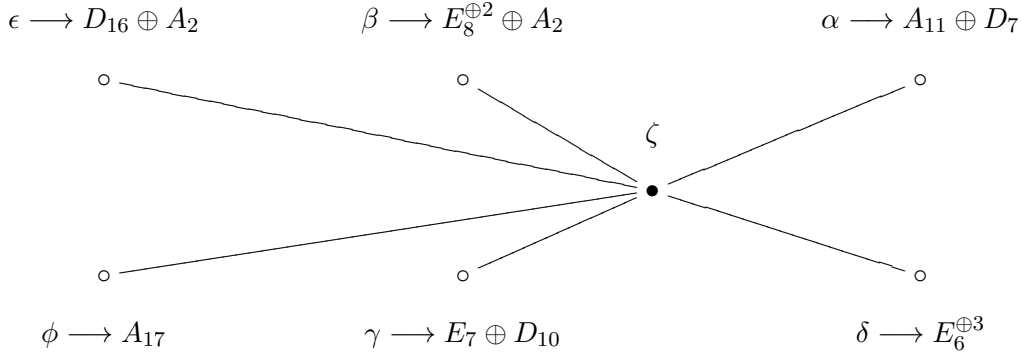
\begin{figure}[htb]
$$\xymatrix@R=.25cm{
{\epsilon \longrightarrow D_{16}\oplus A_2}&&{\beta \longrightarrow E_{8}^{\oplus 2}\oplus A_2}&                      &&{\alpha \longrightarrow A_{11}\oplus D_7}             \\
{\circ}\ar@{-}[ddrrr]&&{\circ}\ar@{-}[ddr]&                      &&{\circ}\ar@{-}[ddll]                \\
&&&                  {\zeta}\\
&&&                {\bullet}      && \\
&  &&                      &&                    \\
{\circ}\ar@{-}[uurrr]&&{\circ}\ar@{-}[uur]&                      &&{\circ}\ar@{-}[uull]                  \\
{\phi \longrightarrow A_{17}}  &&{\gamma \longrightarrow E_7\oplus D_{10}}  &                      &&{\delta \longrightarrow E_6^{\oplus 3}}                     \\
}
$$
\caption{Matching of GIT and Baily--Borel boundary components}\label{bbboundary}
\end{figure}

The fact that the number of Type II  boundary components is equal to the number of GIT boundary components is explained by the fact that a generic point $x_0$ in a GIT boundary component $\alpha$--$\phi$ corresponds to a Type II fourfold $X_0$. By our results (esp. Thm. \ref{mainthm2}) and the general theory, it follows that the period map extends at $x_0$ mapping it to a point  $y_0$ on a Type II boundary component of $(\calD/\Gamma)^*$. This gives a matching between the GIT boundary components and the  Type II Baily--Borel boundary components. To compute the actual matching, we note that  the point $y_0$ encodes the information contained in the graded pieces of the limit mixed Hodge structure $\sH^4_{\lim}$ of a degeneration to $X_0$ (see \cite[Ch. V]{griffithsbook}). This information consists of a discrete part that determines the boundary component and of a continuous part (essentially a $j$-invariant) coming from the Hodge structure on $\Gr_3^W(\sH^4_{\lim})$. In our situation, both pieces of information can be determined in terms of $X_0$ (see section \ref{sectmonodromy}).  We obtain a correspondence between the boundary components in the two compactifications as given in figure \ref{bbboundary} (N.B. $\circ$ and $\bullet$ correspond to  Type II and III components respectively). We omit the details of the computation. Here, we only give the following heuristic argument.  A non-simple singularity  for $X_0$  produces a lattice  of type A-D-E by the following rules:
\begin{itemize}
\item[(1)] an $\widetilde{E}_r$ singularity for $X_0$ gives an $E_r$ lattice;
\item[(2)] an elliptic curve $E$ of degree $d$ in the singular locus of $X_0$ gives an $A_{3d-1}$ lattice  (N.B. $3d$ is the number of special points on $E$, i.e. the points cut on $E$ by  a pencil of cubics degenerating to $X_0$);
\item[(3)] a rational curve $C$ of degree $d$ in the singular locus of $X_0$ gives a $D_{3d+4}$ lattice  (N.B. $3d+4$  special points on $C$).
\end{itemize} 
The direct sum of the lattices associated to the singularities of $X_0$ determines the label of the boundary component associated to $X_0$ (see the discussion for  K3 surfaces from \cite[Remark 5.6]{friedmanannals}).

\subsection{Intersections of hyperplanes from $\calH_{\infty}$} By the work of Looijenga \cite{looijengacompact}, the structure of the birational map $\overline{\calM}\dashrightarrow (\calD/\Gamma)^*$ is very explicit provided that we understand how the hyperplanes from the arrangement $\calH_{\infty}$ intersect inside the period domain $\calD$. This is quite easily achieved by means of the following observation: {\it a hyperplane $H$ from $\calH_{\infty}$ is the restriction to $\calD$ of  a hyperplane $H'\subset \calD'$ orthogonal to a root}, where $\calD'$ is a certain  Type IV domain containing as a subdomain the period domain $\calD$. Specifically, we introduce:

\begin{notation}
Let $\Lambda':=\mathrm{II}_{26,2}$ be the unique even unimodular lattice of signature $(26,2)$. Up to isometries there exists a unique primitive embedding $\Lambda_0\hookrightarrow \Lambda'$. We fix such an embedding and denote by $M\cong E_6$ the orthogonal complement of $\Lambda_0$ in $\Lambda'$. Let  $\calD'$ be the type IV domains associated to $\Lambda'$ and $\calD\hookrightarrow \calD'$ the inclusion induced by the embedding $\Lambda_0\hookrightarrow \Lambda'$. 
\end{notation}

With these notations, we have:

\begin{lemma}\label{lemmaroot}
Let $H$ be a hyperplane from the arrangement $\calH_{\infty}$. Then there exists a root $\delta'\in \Lambda'$ (i.e. $\delta'^2=2$) such that
\begin{itemize}
\item[i)] the sublattice in $\Lambda'$ spanned by $M$ and $\delta'$  is isomorphic to $E_7$; 
\item[ii)] the restriction to $\calD$ of the hyperplane $H'_{\delta'}\subset \calD'$  orthogonal  to $\delta'$  coincides with $H$.
\end{itemize} 
Conversely, for any root $\delta'\in \Lambda'$ satisfying the condition i), the restriction to $\calD$ of the associated orthogonal hyperplane $H'_{\delta'}$ is an element of $\calH_{\infty}$.
\end{lemma}
\begin{proof}
As noted in section \ref{sectperiod}, given $H$ in $\calH_{\infty}$, there exists a long root $\delta\in \Lambda_0$ orthogonal to $H$. The condition of long root implies in particular that $\frac{\delta}{3}$ is a generator of the discriminant group $A_{\Lambda_0}$. Since $M$ and $\Lambda_0$ are glued along the discriminant $A_{\Lambda_0}$ to give the lattice $\Lambda'$, it follows easily that the saturation $S$ of the sublattice of $\Lambda'$ spanned by $M\cong E_6$ and $\delta$ is isomorphic to $E_7$. Clearly, the orthogonal hyperplane $H'_{\delta'}$ to a root  $\delta'\in S\setminus M$  restricts on $\calD$ to $H$. 
Since there is a unique embedding of $E_6$ into $E_7$ and the orthogonal complement of $E_6$ in this embedding is the rank $1$ lattice spanned by an element of norm $6$, the converse also follows.
\end{proof}

A similar description  holds also for the hyperplanes from $\calH_{\Delta}$. The only difference is to replace $E_7$ from condition i) by $E_6\oplus A_1$. In particular, from the work of Borcherds, it follows that $(\calH_{\infty}\cup \calH_{\Delta})/\Gamma$ is the zero locus of an automorphic form. Thus, $(\calD\setminus (\calH_{\infty}\cup \calH_{\Delta}))/\Gamma$ is an affine variety. This follows also from GIT via the isomorphism $(\calD\setminus (\calH_{\infty}\cup \calH_{\Delta}))/\Gamma\cong \calM_0$.

\begin{proposition}
There exists an automorphic form $\Phi$ on $\calD$ such that 
\begin{itemize}
\item[i)] the weight of $\Phi$ is $48$;
\item[ii)] the vanishing locus of $\Phi$ is $\calH_{\infty}\cup \calH_{\Delta}$;
\item[iii)] $\phi$ vanishes of order $1$ along $\calH_{\Delta}$ and of order $27$ along $\calH_{\infty}$.
\end{itemize}
We call $\Phi$ the discriminant automorphic form.
\end{proposition}
\begin{proof}
The situation is formally the same as for degree two K3 surface (see Borcherds et al. \cite[Thm. 1.2, Ex. 2.1]{bkpsb}). Namely, it is well known that  on $\calD'$ there exists an automorphic form $\Phi_{12}$  vanishing exactly along the hyperplanes $H'_{\delta'}$, where $\delta'$ is a root of $\Lambda'$. By restricting $\Phi_{12}$ to $\calD$ and dividing by the product of the linear forms which vanish along the hyperplanes $H'_{\delta'}$ containing $\calD$,  we obtain an automorphic form $\Phi$ on $\calD$. The weight of $\Phi$ is the weight of $\Phi_{12}$ plus  the number of hyperplanes $H'_{\delta'}$ containing $\calD$ (i.e. half the number of roots of $E_6$). Clearly, the vanishing locus of $\Phi$ is the union of hyperplanes $H$ which are the restrictions to $\calD$ of the hyperplanes of type $H'_{\delta'}$ for $\delta'$ a root of $\Lambda'$. The condition that $H'_{\delta'}\cap \calD\neq \emptyset$ implies that the sublattice $\langle M, \delta'\rangle_{\Lambda'}$ (the span of $M$ and $\delta'$) is positive definite. Since $M\cong E_6$, we obtain that  $\langle M, \delta'\rangle_{\Lambda'}$ is isomorphic to either  $E_6\oplus A_1$ or $E_7$. The two cases cases correspond to  $\calH_{\Delta}$ and $\calH_{\infty}$ respectively. The claim follows.
\end{proof}

\begin{remark}
The advantage of working with $\Lambda'$ instead of $\Lambda_0$ or $\Lambda$ is that the similarities and differences between the arrangements $\calH_{\Delta}$ and $\calH_{\infty}$ are more apparent. Also, since we are working with an even unimodular lattice, the results of Nikulin \cite{nikulin} are easier to apply. 
\end{remark}
An immediate consequence of lemma \ref{lemmaroot} is that a non-empty intersection in $\calD$ of  hyperplanes from $\calH_\infty$ has codimension at most $2$ and that all codimension $2$ intersections are conjugated by $\Gamma$. 
\begin{lemma}\label{intersect1}
Assume that $H_1,\dots, H_k$ (with $k\ge 2$) are distinct hyperplanes from $\calH_{\infty}$ such that $\calD\cap H_1\cap\dots\cap H_k\neq \emptyset$, then the locus $\calD\cap H_1\cap\dots\cap H_k$ has codimension $2$ in $\calD$. Furthermore, all codimension two intersection loci of hyperplanes from $\calH_{\infty}$ are conjugate by $\Gamma$.
\end{lemma}
\begin{proof}
By lemma \ref{lemmaroot}, each hyperplane $H_i$ corresponds to a root $\delta_i'\in \Lambda'$. It follows that the sublattice of $\Lambda'$ spanned by $M$ and the roots $\delta_i'$ is a root lattice $R$ of rank at least $8$ (the hyperplanes $H_i$ are distinct). The condition that  $\calD\cap H_1\cap\dots\cap H_k\neq \emptyset$ is equivalent to $R$ being positive definite. Since $\langle M, \delta_i'\rangle_{\Lambda'}\cong E_7$ for each $i$, it follows that $R=\langle M,\delta_1',\dots, \delta_k'\rangle_{\Lambda'}$ is an irreducible positive definite root lattice of rank at least $8$. Since $M\cong E_6$ embeds in $R$, we conclude that $R\cong E_8$. The codimension claim follows. Since the embeddings $E_6\hookrightarrow E_8$, $E_8\hookrightarrow \Lambda'$, and  $E_6\hookrightarrow \Lambda'$ are unique up to isometries, it is not hard to see that the codimension $2$ intersections form a single $\Gamma$-orbit.
\end{proof}
 
Geometrically, the previous lemma implies that the period map for cubic fourfolds is undefined precisely along the curve $\chi$. Moreover, via the diagram (\ref{birationaldiag}), the curve $\chi$ is a flip of the locus in $\calD/\Gamma$ corresponding to the codimension $2$ intersections from $\calH_{\infty}$ and the point $\omega$ is a contraction of the divisor $\calH_{\infty}/\Gamma$.
 
 \smallskip
 
To understand the effect of the small birational modification $\widehat{\calM}\to (\calD/\Gamma)^*$, we have to understand the structure of the arrangement $\calH_{\infty}$  near the boundary components. Since a Baily--Borel boundary component corresponds to an isotropic sublattice $E$ of $\Lambda_0$, this is essentially equivalent to understanding how the hyperplanes from $\calH_{\infty}$ containing $E$ intersect. Specifically, one is interested in the following lemma:
\begin{lemma}\label{intersect2}
Let $E$ be a rank $2$ isotropic primitive sublattice of $\Lambda_0$ and $R$ the root sublattice in $E^\perp/E$ (see Prop. \ref{thmbbcompact}). Denote by $K$ the restriction of the intersection of hyperplanes from  $\calH_\infty$ containing $E$ to $E^\perp$,  i.e.
$$K=\cap_{H\in \calH_{\infty},\ E\subset H}\ (H\cap E^\perp)$$ 
Then, the following holds:
\begin{itemize}
\item[i)]($E_6$ type) if $R$ is either $E_6^{\oplus 3}$ or $A_{11}\oplus D_7$, then $K=E^\perp$;
\item[ii)]($E_7$ type)  if $R$ is either $E_7\oplus D_{10}$ or $A_{17}$, then $K$ has codimension $1$ in $E^\perp$;
\item[iii)]($E_8$ type)  if $R$ is either $E_8^{\oplus 2}\oplus A_2$ or $D_{16}\oplus A_2$, then $K$ has codimension $2$ in $E^\perp$.
\end{itemize}
\end{lemma}
\begin{proof} As before, a hyperplane $H$ from $\calH_\infty$ corresponds to a root $\delta'$ in $\Lambda'$ such that $\langle M, \delta'\rangle_\Lambda'\cong E_7$. The condition $E\subset H$ is equivalent to $\delta'$ being orthogonal to $E$. By the proof of \ref{thmbbcompact}, it follows that  $L=E^\perp_{\Lambda'}/E$ is an even unimodular lattice of rank $24$ whose isometry class is uniquely determined by $E$.  Additionally, the projection of $M\subset E^\perp_{\Lambda'}$ in $L$ determines an embedding $E_6\hookrightarrow L$ (uniquely determined by $E$). Similarly, $\delta'$ projects to a root of $L$. Now, the claim follows by  arguments similar to those of \ref{intersect1} (with $L$ in place of $\Lambda'$). \end{proof}

 The role of linear space $K$ associated to a boundary component in the construction of the blow-up  $\widehat{\calM}\to (\calD/\Gamma)^*$ is explained in  Looijenga \cite[\S3]{looijengacompact}. In particular, the following geometric consequences follow immediately from the general theory and the computations of the lemma:
\begin{itemize}
\item[i)] The GIT boundary components $\alpha$ and $\delta$ are unaffected by the birational transformation from (\ref{birationaldiag}). In particular, these strata are $1$-dimensional and the period map extends to an  isomorphism over these strata.
\item[ii)] The strata $\gamma$ and $\phi$  are $2$-dimensional and the only indeterminacy point for the period map on these strata is the point $\omega$.
\item[iii)] The strata $\beta$ and $\epsilon$ are $3$-dimensional, and the indeterminacy locus along these strata is the curve $\chi$ (including the special point $\omega$). 
\end{itemize}
Additionally, the open subsets in $\alpha$--$\phi$ parametrizing Type II cubic fourfolds (i.e. excluding the surface $\sigma\subset \alpha\cup \dots \cup\phi$ of  \cite{gitcubic}) fiber over the Type II boundary components of the Baily--Borel compactification (see  the matching of components from figure \ref{bbboundary}). Clearly, the dimensions of the fibers are $2$ for $\beta$ and $\epsilon$, $1$ for $\gamma$ and $\phi$, and $0$ for  $\alpha$ and $\delta$ respectively. Each of the Type II boundary components in $(\calD/\Gamma)^*$ are quotients of the upper half space by some arithmetic group. Thus, the fibration mentioned above is obtained (up to some finite index)  geometrically by associating a natural $j$-invariant to a Type II fourfold. For example, the Type II fourfolds parameterized by the boundary component $\phi$ are singular along an elliptic normal curve of degree $6$. The $j$-invariant mentioned above is the $j$-invariant of this elliptic curve. A similar situation holds also for the other boundary components: there is a natural $j$-invariant associated to the singularities of a Type II fourfold (see Prop. \ref{propsingtype} (i) and Remark \ref{remresolution}).

%%%%%%%%%%%%%%%%%%%%%%%%%%%%%%%%%%%%%%%%%%%%%%%%%%%%%%%%%%%%%%
\bibliography{references}

\def\cprime{$'$}
\providecommand{\bysame}{\leavevmode\hbox to3em{\hrulefill}\thinspace}
\providecommand{\MR}{\relax\ifhmode\unskip\space\fi MR }
\providecommand{\MRhref}[2]{%
  \href{http://www.ams.org/mathscinet-getitem?mr=#1}{#2}
}
\providecommand{\href}[2]{#2}
\begin{thebibliography}{10}

\bibitem{allcock1}
D.~Allcock, \emph{The moduli space of cubic threefolds}, J. Algebraic Geom.
  \textbf{12} (2003), no.~2, 201--223.

\bibitem{allcock3fold}
D.~Allcock, J.~A. Carlson, and D.~Toledo, \emph{{The Moduli Space of Cubic
  Threefolds as a Ball Quotient}},  Mem. Amer. Math. Soc. 209(985):xii+70, 2011.
  
\bibitem{AGV1}
V.~I. Arnol{\cprime}d, S.~M. Guse{\u\i}n-Zade, and A.~N. Varchenko,
  \emph{Singularities of differentiable maps. {V}ol. {I}}, Monographs in
  Mathematics, vol.~82, Birkh\"auser, Boston, MA, 1985.

\bibitem{beauvilleprym}
A.~Beauville, \emph{Vari\'et\'es de {P}rym et jacobiennes interm\'ediaires},
  Ann. Sci. \'Ecole Norm. Sup. (4) \textbf{10} (1977), no.~3, 309--391.

\bibitem{beauvillemon}
\bysame, \emph{Le groupe de monodromie des familles universelles
  d'hypersurfaces et d'intersections compl\`etes}, Complex analysis and
  algebraic geometry (G\"ottingen, 1985), Lecture Notes in Math., vol. 1194,
  Springer, Berlin, 1986, pp.~8--18.

\bibitem{beauvilledonagi}
A.~Beauville and R.~Donagi, \emph{La vari\'et\'e des droites d'une hypersurface
  cubique de dimension {$4$}}, C. R. Acad. Sci. Paris S\'er. I Math.
  \textbf{301} (1985), no.~14, 703--706.

\bibitem{bkpsb}
R.~E. Borcherds, L.~Katzarkov, T.~Pantev, and N.~I. Shepherd-Barron,
  \emph{Families of {$K3$} surfaces}, J. Algebraic Geom. \textbf{7} (1998),
  no.~1, 183--193.

\bibitem{carlsonperiod}
J.~A. Carlson, S.~M{\"u}ller-Stach, and C.~A.~M. Peters, \emph{Period mappings
  and period domains}, Cambridge Studies in Advanced Mathematics, vol.~85,
  Cambridge University Press, Cambridge, 2003.

\bibitem{conway}
J.~H. Conway and N.~J.~A. Sloane, \emph{Sphere packings, lattices and groups},
  third ed., Grundlehren der Mathematischen Wissenschaften, vol. 290,
  Springer-Verlag, New York, 1999.

\bibitem{sga72}
P.~Deligne and N.~Katz (eds.), \emph{Groupes de monodromie en g\'eom\'etrie
  alg\'ebrique. {II}}, Springer-Verlag, Berlin, 1973, S\'eminaire de
  G\'eom\'etrie Alg\'ebrique du Bois-Marie 1967--1969 (SGA 7 II), Lecture Notes
  in Math., Vol. 340.

\bibitem{dimca}
A.~Dimca, \emph{Sheaves in topology}, Universitext, Springer-Verlag, Berlin,
  2004.

\bibitem{dolgachevinsignificant}
I.~V. Dolgachev, \emph{Cohomologically insignificant degenerations of algebraic
  varieties}, Compositio Math. \textbf{42} (1980/81), no.~3, 279--313.

\bibitem{versality}
A.~A. du~Plessis, \emph{Versality properties of projective hypersurfaces}, Real
  and Complex Singularities, Trends Math., Birkh\"auser, Basel, 2006,
  pp.~289--298.

\bibitem{duplessiswall0}
A.~A. du~Plessis and C.~T.~C. Wall, \emph{Singular hypersurfaces, versality,
  and {G}orenstein algebras}, J. Algebraic Geom. \textbf{9} (2000), no.~2,
  309--322.

\bibitem{einsb}
L.~Ein and N.~I. Shepherd-Barron, \emph{Some special {C}remona
  transformations}, Amer. J. Math. \textbf{111} (1989), no.~5, 783--800.

\bibitem{friedmanannals}
R.~Friedman, \emph{A new proof of the global {T}orelli theorem for {$K3$}
  surfaces}, Ann. of Math. (2) \textbf{120} (1984), no.~2, 237--269.

\bibitem{fultonintersection}
W.~Fulton, \emph{Intersection theory}, second ed., Ergebnisse der Mathematik
  und ihrer Grenzgebiete. 3. Folge., Springer-Verlag, Berlin, 1998.

\bibitem{fultonharris}
W.~Fulton and J.~Harris, \emph{Representation theory}, Graduate Texts in
  Mathematics, vol. 129, Springer-Verlag, New York, 1991.

\bibitem{griffithsbook}
P.~A. Griffiths (ed.), \emph{Topics in transcendental algebraic geometry},
  Annals of Mathematics Studies, vol. 106, Princeton, NJ, Princeton University
  Press, 1984.

\bibitem{griffithscurvature}
P.~A. Griffiths and L.~Tu, \emph{Curvature properties of the {H}odge bundles},
  Topics in transcendental algebraic geometry, Ann. of Math. Stud., vol. 106,
  Princeton Univ. Press, 1984, pp.~29--49.

\bibitem{loeser}
G.~Guibert, F.~Loeser, and M.~Merle, \emph{Iterated vanishing cycles,
  convolution, and a motivic analogue of a conjecture of {S}teenbrink}, Duke
  Math. J. \textbf{132} (2006), no.~3, 409--457.

\bibitem{hassett}
B.~Hassett, \emph{Special cubic fourfolds}, Compositio Math. \textbf{120}
  (2000), no.~1, 1--23, (long version available at
  http://math.rice.edu/~hassett/papers.html).

\bibitem{huybrechts}
D.~Huybrechts, \emph{Compact hyper-{K}\"ahler manifolds: basic results},
  Invent. Math. \textbf{135} (1999), no.~1, 63--113.

\bibitem{kirwanhyp}
F.~Kirwan, \emph{Moduli spaces of degree {$d$} hypersurfaces in {${\bf P}\sb
  n$}}, Duke Math. J. \textbf{58} (1989), no.~1, 39--78.

\bibitem{kirwan}
F.~C. Kirwan, \emph{Partial desingularisations of quotients of nonsingular
  varieties and their {B}etti numbers}, Ann. of Math. (2) \textbf{122} (1985),
  no.~1, 41--85.

\bibitem{ksb}
J.~Koll{\'a}r and N.~I. Shepherd-Barron, \emph{Threefolds and deformations of
  surface singularities}, Invent. Math. \textbf{91} (1988), no.~2, 299--338.

\bibitem{kulikovmhs}
V.~S. Kulikov, \emph{Mixed {H}odge structures and singularities}, Cambridge
  Tracts in Mathematics, vol. 132, Cambridge University Press, Cambridge, 1998.

\bibitem{gitcubic}
R.~Laza, \emph{The moduli space of cubic fourfolds}, J. Algebraic Geom.
  \textbf{18} (2009), no.~3, 511--545.

\bibitem{looijengatriangle}
E.~Looijenga, \emph{The smoothing components of a triangle singularity. {II}},
  Math. Ann. \textbf{269} (1984), no.~3, 357--387.

\bibitem{looijengacompact}
\bysame, \emph{Compactifications defined by arrangements. {II}. {L}ocally
  symmetric varieties of type {IV}}, Duke Math. J. \textbf{119} (2003), no.~3,
  527--588.

\bibitem{looijengacubic}
\bysame, \emph{{The period map for cubic fourfolds}}, Invent. Math. \textbf{177}
  (2009), no.~1, 213--233.
  
  \bibitem{looijengaswierstra}
E.~Looijenga and R.~Swierstra, \emph{{The period map for cubic threefolds}}, Compositio Math.
  \textbf{143} (2007), no.~4, 1037--1049.
  
\bibitem{looijengaswierstra1}
 \bysame, \emph{{On period maps that are open
  embeddings}}, J. Reine Angew. Math., \textbf{617} (2008), 169--192.



\bibitem{clemensschmid}
D.~R. Morrison, \emph{The {C}lemens-{S}chmid exact sequence and applications},
  Topics in transcendental algebraic geometry, Ann. of Math. Stud., vol. 106,
  Princeton Univ. Press, 1984, pp.~101--119.

\bibitem{mumford}
D.~Mumford, \emph{Stability of projective varieties}, Enseignement Math. (2)
  \textbf{23} (1977), no.~1-2, 39--110.

\bibitem{GIT}
D.~Mumford, J.~Fogarty, and F.~Kirwan, \emph{Geometric invariant theory}, third
  ed., Ergebnisse der Mathematik und ihrer Grenzgebiete (2), vol.~34,
  Springer-Verlag, Berlin, 1994.

\bibitem{nikulin}
V.~V. Nikulin, \emph{Integral symmetric bilinear forms and some of their
  applications}, Math. USSR Izvestiya \textbf{43} (1980), no.~1, 103--167.

\bibitem{ogrady}
K.~G. O'Grady, \emph{Irreducible symplectic 4-folds numerically equivalent to
  {$(K3)^{[2]}$}}, arXiv:math/0504434 [math.AG] (2005), 44 pp.

\bibitem{peterssteenbrink}
C.~Peters and J.~H.~M. Steenbrink, \emph{{M}ixed {H}odge {S}tructures},
Ergebnisse der
  Mathematik und ihrer Grenzgebiete (3),  vol.~52, Springer-Verlag, Berlin, 2008.

\bibitem{pinkhamduality}
H.~C. Pinkham, \emph{Singularit\'es exceptionnelles, la dualit\'e \'etrange
  d'{A}rnold et les surfaces {$K-3$}}, C. R. Acad. Sci. Paris S\'er. A-B
  \textbf{284} (1977), no.~11, A615--A618.

\bibitem{saitosc}
M.~Saito, \emph{On {S}teenbrink's conjecture}, Math. Ann. \textbf{289} (1991),
  no.~4, 703--716.

\bibitem{scattone}
F.~Scattone, \emph{On the compactification of moduli spaces for algebraic
  {$K3$} surfaces}, Mem. Amer. Math. Soc. \textbf{70} (1987), no.~374, x+86.

\bibitem{shahinsignificant}
J.~Shah, \emph{Insignificant limit singularities of surfaces and their mixed
  {H}odge structure}, Ann. of Math. (2) \textbf{109} (1979), no.~3, 497--536.

\bibitem{shah}
\bysame, \emph{A complete moduli space for {$K3$} surfaces of degree {$2$}},
  Ann. of Math. (2) \textbf{112} (1980), no.~3, 485--510.

\bibitem{siersma}
D.~Siersma, \emph{The vanishing topology of non isolated singularities}, New
  developments in singularity theory (Cambridge, 2000), vol.~21, Kluwer Acad.
  Publ., 2001, pp.~447--472.

\bibitem{steenbrinkmhs}
J.~H.~M. Steenbrink, \emph{Mixed {H}odge structure on the vanishing
  cohomology}, Real and complex singularities (Oslo, 1976), Sijthoff and
  Noordhoff, Alphen aan den Rijn, 1977, pp.~525--563.

\bibitem{steenbrinkinsignificant}
\bysame, \emph{Cohomologically insignificant degenerations}, Compositio Math.
  \textbf{42} (1980/81), no.~3, 315--320.

\bibitem{steenbrinksemicontinuity}
\bysame, \emph{Semicontinuity of the singularity spectrum}, Invent. Math.
  \textbf{79} (1985), no.~3, 557--565.

\bibitem{voisin}
C.~Voisin, \emph{Th\'eor\`eme de {T}orelli pour les cubiques de {${\bf P}\sp
  5$}}, Invent. Math. \textbf{86} (1986), no.~3, 577--601.

\bibitem{wallsextic}
C.~T.~C. Wall, \emph{Sextic curves and quartic surfaces with higher
  singularities}, preprint (1998), 31 pp.

\end{thebibliography}
\end{document}